\documentclass[sn-mathphys,autheryear]{sn-jnl}

\usepackage{multirow}
\usepackage{subfigure}
\usepackage{amsmath,amsfonts,amssymb,amsthm,booktabs,color,epsfig,graphicx,hyperref,url}
\usepackage{enumerate}
\usepackage{lineno,hyperref}
\usepackage{braket}
\usepackage{enumitem}
\usepackage{tabularcalc}

\numberwithin{equation}{section}

\theoremstyle{plain}
\newtheorem{thm}{Theorem}

\newtheorem{rmk}{Remark}
\newtheorem{prop}{Proposition}

\newtheorem{hypo}{Hypothesis}

\theoremstyle{definition}

\usepackage{graphicx}%
\usepackage{multirow}%
\usepackage{amsmath,amssymb,amsfonts}%
\usepackage{amsthm}%
\usepackage{mathrsfs}%
\usepackage[title]{appendix}%
\usepackage{xcolor}%
\usepackage{textcomp}%
\usepackage{manyfoot}%
\usepackage{booktabs}%
\usepackage{algorithm}%
\usepackage{algorithmicx}%
\usepackage{algpseudocode}%
\usepackage{listings}%
\usepackage{color}
\usepackage {lineno}

\theoremstyle{thmstyleone}%
\theoremstyle{thmstyletwo}%

\theoremstyle{thmstylethree}%

\begin{document}

\switchlinenumbers
\allowdisplaybreaks

\title[]{Statistical inference for multi-regime threshold Ornstein$\mathit{-}$Uhlenbeck processes}


\author[1]{\fnm{Yuecai} \sur{Han}}\email{hanyc@jlu.edu.cn}
\equalcont{These authors contributed equally to this work.}

\author*[1]{\fnm{Dingwen} \sur{Zhang}}\email{zhangdw20@mails.jlu.edu.cn}
\equalcont{These authors contributed equally to this work.}

\affil[1]{\orgdiv{School of Mathematics}, \orgname{Jilin University}, \orgaddress{\street{Qianjin Street}, \city{Changchun}, \postcode{130012}, \state{Jilin}, \country{China}}}

\abstract{In this paper, we investigate the parameter estimation for threshold Ornstein$\mathit{-}$Uhlenbeck processes. Least squares method is used to obtain continuous-type and discrete-type estimators for the drift parameters based on continuous and discrete observations, respectively. The strong consistency and asymptotic normality of the proposed least squares estimators are studied. We also propose a modified quadratic variation estimator based on the long-time observations for the diffusion parameters and prove its consistency. { Our simulation results suggest that the performance of our proposed estimators for the drift parameters may show improvements compared to generalized moment estimators. Additionally, the proposed modified quadratic variation estimator exhibits potential advantages over the usual quadratic variation estimator with relatively small sample sizes.}
	In particular, our method can be applied to the multi-regime cases ($m>2$), while the generalized moment method only deals with the two regime cases ($m=2$). The U.S. treasury rate data is used to illustrate the theoretical results.}




\maketitle
 
\section{Introduction}\label{sec1}
Let $(\Omega, \mathcal{F}, \mathbb{P})$ be a complete probability space equipped with a right continuous and increasing family of $\sigma$-algebra $\{\mathcal{F}_{t}\}_{t\geq0}$ and $W=\{W_{t}\}_{t\geq0}$ be a one-dimensional Brownian motion adapted to $\{\mathcal{F}_{t}\}_{t\geq0}$. The threshold Ornstein$\mathit{-}$Uhlenbeck (OU) process $X=\{X_{t}\}_{t\geq0}$ is defined as the unique solution to the following threshold stochastic differential equation (SDE)
\begin{equation}\label{eq1}
	\left\{
	\begin{aligned}
		&\mathrm{d}X_{t}=\sum_{j=1}^{m}(\beta_{j}-\alpha_{j}X_{t})I_{\{\theta_{j-1}\leq X_{t}<\theta_{j}\}}\mathrm{d}t+\sigma \mathrm{d}W_{t},\\
		&X_{0}=x,
	\end{aligned}
	\right.
\end{equation}
where $\beta_{j}, \alpha_{j}\in\mathbb{R}$ are the drift parameters, $\Theta=(\theta_{0}, \cdots, \theta_{m})$ satisfy $-\infty=\theta_{0}<\theta_{1}<\cdots<\theta_{m}=\infty$ are the so-called thresholds, $\sigma>0$ is the diffusion parameter, $x\in\mathbb{R}$ is the initial value of the process $\{X_{t}\}_{t\geq0}$, and $I_{\{\cdot\}}$ denotes the indicator function. Moreover, we assume that $\alpha_{1}$ and $\alpha_{m}$ are larger than zero to ensure the ergodicity of $\{X_{t}\}_{t\geq0}$. The existence and uniqueness of the solution to Eq. (\ref{eq1}) follow, for example, from the results of \cite{Bass and Pardoux(1987)}. 
Suppose that the discretely observed process is observed at the regularly spaced time points $\{t_{k}=kh\}_{k=0}^{n}$ where $h$ is the mesh size and $n$ is the sample size. In this paper, we aim to propose the least squares estimators (LSEs) for the unknown parameters $\alpha_{j}$ and $\beta_{j}$ based on either the continuously observed process $\{X_{t}\}_{t\geq0}$ or discretely observed process $\{X_{t_{kh}}\}_{k=0}^{n}$, respectively. Furthermore, we propose the quadratic variation estimator (QVE) and the modified quadratic variation estimator (MQVE) for the diffusion parameters.

As useful stochastic dynamics models, threshold diffusion processes have been extensively employed in numerous fields. The threshold autoregressive models are introduced to model the nonlinearities in nonlinear time series \citep[see, e.g.,][]{Brockwell and Hyndman(1991),Brockwell and Hyndman(1992)}. In finance, \cite{Lejay and Pigato(2019)} show that an exponential form of the threshold process generalizes the Black-Scholes model in a way to model leverage effects. \cite{Gerber and Shiu(2006)} model the surplus of a company after the payment of dividends, which are paid only if the profits of the company are higher than a certain threshold. The threshold processes also have numerous applications in physics \citep[see, e.g.][]{Ramirez et al.(2013),Sattin(2008)}, meteorology \citep[see, e.g.][]{Hottovy and Stechmann(2015)}, forecast \citep[see, e.g.][]{Brockwell and Hyndman(1992)}, and they are in close contact with the skew diffusion processes \citep[see, e.g.][]{Ding et al.(2021),Gairat and Shcherbakov(2017)}.

For threshold OU processes, \cite{Chan(1993)} proposes a LSE for a stationary ergodic threshold autoregressive model and proves its consistency and limiting distribution. \cite{Su and Chan(2015),Su and Chan(2017)} study the asymptotic behavior of the quasi-likelihood estimator of a diffusion with piecewise regular diffusivity and piecewise affine drift with an unknown threshold. They provide a hypothesis test to decide whether or not the drift is affine or piecewise affine. \cite{Lejay and Pigato(2020)} propose the maximum likehood estimator (MLE) for the drift parameters, with both drift and diffusion coefficients constant on the positive and negative axis, yet discontinuous at zero.  \cite{Kutoyants(2012)} studies the asymptotic behavior of the maximum likelihood estimator and Bayesian estimator for the threshold parameter. He also discusses the possibility of the construction of the goodness-of-fit test. \cite{Mazzonetto and Pigato(2020)} discuss the (quasi)-MLE of the drift parameters, both assuming continuous and discrete time observations. \cite{Hu and Xi(2022)} propose a generalized moment estimator (GME) to estimate the drift parameter based on the discretely observed process. Their main theoretical basis is the celebrated ergodic theorem. \cite{Han and Zhang(2023)} propose the trajectory fitting estimators for the drift parameters and obtain the asymptotic behavior of the estimators. We also draw attention to a related piece of work by  \cite{Xing et al.(2020)}, in which they propose an estimator for the drift parameter of a skew OU process.

{Based on the celebrated Girsanov's theorem for semimartingales, one could obtain the log-likehood function, and the MLE for the drift parameter is derived by minimizing the log-likehood function \citep{Su and Chan(2015),Su and Chan(2017),Mazzonetto and Pigato(2020)}.} The generalized moment method could deal with low-frequency data such as the gross domestic product. However, it is significantly more complicated to explore the case $m>2$ \cite[see][]{Hu and Xi(2022)}. In this paper, we apply the least squares method to deal with the drift parameter estimation problem based on a continuously observed process and a discretely observed process, respectively. 

The remainder of this paper is organized as follows. In Section \ref{secpre}, we introduce some basic facts on the stationary and ergodic properties of threshold OU processes. In Section \ref{secc} and \ref{secd}, we construct the LSEs for the drift parameters based on continuous sampling and discrete sampling. We prove the strong consistency and asymptotic normality of the proposed LSEs. In section \ref{secsigma}, we consider a general threshold OU process, and propose the QVEs and the MQVEs for the diffusion parameters $\sigma$. The results of our simulation studies are reported in Section \ref{secnr}, exhibiting improved performance versus the generalized moment method in two-regime cases. In the case of multi-regime, our proposed LSEs and MQVEs have great performance. We also show an application of our method in the U.S. treasury rate. Section \ref{seccon} concludes with some discussion on the further work.

\section{Preliminaries}\label{secpre}
In this section, we first introduce some basic facts on the stationary and ergodic properties of threshold OU processes. 
Throughout this paper, we shall
use the notation ``$\stackrel{a.s.}{\longrightarrow}$'' to denote ``almost sure convergence'', the notation ``$\stackrel{P}{\longrightarrow}$" to denote ``convergence in probability'', {and the notation ``$X_{T}\sim \mathcal{N}(a,b)$'' or ``$X_{n}\sim \mathcal{N}(a,b)$'' to denote ``the convergence in distribution of $X_{T}$ or $X_{n}$ to $X$ where $X$ follows a normal distribution with mean $a$ and variance $b$''.} For the convenience of the following description, we make a notation $I_{\{\cdot,j\}}=I_{\{\theta_{j-1}\leq \cdot<\theta_{j}\}}$, $j=1,\cdots, m$. 

The following proposition, adopt from Proposition 3.1 of \cite{Brockwell and Hyndman(1991)}, Section 1 of \cite{Browne and Whitt(1995)}, and Proposition 1 of \cite{Hu and Xi(2022)}, shows that the process $\{X_{t}\}_{t\geq0}$ has a unique invariant density. 

\hfill

\begin{prop}\label{prop1}
	If the drift parameters $\alpha_{1}$, $\alpha_{m}$, $\beta_{1}$, and $\beta_{m}$ satisfy
	\begin{equation*}
		\lim _{x \rightarrow-\infty}\left(-\alpha_{1} x^{2}+2 \beta_{1} x\right)<0\quad\text{and}\quad \quad \lim _{x \rightarrow \infty}\left(-\alpha_{m} x^{2}+2 \beta_{m} x\right)<0,
	\end{equation*}
	the unique invariant density of the process $\{X_{t}\}_{t\geq0}$  is given by
	\begin{equation}\label{equid1}
		\pi(x)=\sum_{j=1}^{m}\pi_{j}(x)=\sum_{j=1}^{m} k_{j} \exp \left(\left(-\alpha_{j} x^{2}+2 \beta_{j} x\right) / \sigma^{2}\right) I_{\{x,j\}},
	\end{equation}
	where $k_{j}$ are uniquely determined by the following equations
	\begin{equation*}\label{equid2}
		\int_{-\infty}^{\infty} \pi(x) d x=1 \quad \text { and } \quad \pi\left(\theta_{j}-\right)=\pi\left(\theta_{j}+\right), \quad j=1,2, \ldots, m-1.
	\end{equation*}
\end{prop}

The following proposition, based on basic stability theories of Markov processes, describes the ergodic properties of threshold OU processes \citep[see Lemma 1 of][]{Hu and Xi(2022)}. 

\hfill

\begin{prop}
	{Assume that the drift parameters $\alpha_{1}$, $\alpha_{m}$, $\beta_{1}$, and $\beta_{m}$satisfy the conditions of Proposition \ref{prop1}.} Then the discretely observed process $\{X_{kh}\}_{k=0}^{n}$ and continuously observed process $\{X_{t}\}_{t\geq0}$ are ergodic, i.e., for any $f\in L_{1}(\mathbb{R})$,
	\begin{equation}\label{eqergodic}
		\lim_{n\rightarrow\infty}\frac{1}{n}\sum_{k=0}^{n}f(X_{t_{kh}})=\lim_{T\rightarrow\infty}\frac{1}{T}\int_{0}^{T}f(X_{t})\mathrm{d}t=\mathbb{E}[f(X_{\infty})]=\int_{-\infty}^{\infty}f(x)\pi(x)\mathrm{d}x, \quad a.s.,
	\end{equation}
	where $\pi(x)$ is defined in Eq. (\ref{equid1}).
\end{prop}

\hfill

\begin{rmk}
	When $\beta_{j}=0$, $j=1,2,\cdots,m$, the sufficient and necessary condition for the ergodic properties of the threshold OU process is 
	\begin{equation*}
		\alpha_{1}>0 \quad\text{and}\quad \alpha_{m}>0.
	\end{equation*}
\end{rmk}

Let $\braket{X,Y}$ be the covariation process of $X$ and $Y$, and $\braket{X}=\braket{X,X}$. The following proposition, consider the representation of martingale, i.e., we treat the continuous local martingale as a time-changed Brownian motion \citep[see Theorem 1.6, Chapter 5, p.188 of][]{Revuz and Yor(1998)}.

\hfill

	\begin{prop}\label{propset}
		If $M$ is a $\mathcal{F}_{t}$-continuous local martingale vanishing at $0$ and such that $\braket{M}_{\infty}=\infty$, then we have the following results.
		\begin{enumerate}[label*=(\arabic*)]
			\item $B_{t}=M_{\tau_{t}}$ is the so-called Dambis, Dubins-Schwarz Brownian motion and $M_{t}=B_{\braket{M}_{t}}$, where
			\begin{equation*}
				\tau_{t}=\inf\{s:\braket{M}_{s}>t\}.
			\end{equation*}
			\item The law of large numbers for continuous local martingale:
			\begin{equation*}
				\lim_{T\rightarrow\infty}\frac{M_{T}}{\braket{M}_{T}}=0,\quad a.s.
			\end{equation*}
			\item The central limit theorem for continuous local martingale:
			\begin{equation*}
				\frac{M_{T}}{\sqrt{\braket{M}_{T}}}\sim\mathcal{N}(0,1),\quad\text{as}\quad T\rightarrow\infty.
			\end{equation*}
			\item The law of the iterated logarithm:
			\begin{equation*}
				\begin{aligned}
					&\varlimsup_{T\rightarrow\infty}\frac{M_{T}}{\sqrt{2\braket{M}_{T}\log\log\braket{M}_{T}}}=1,\\
					&\varliminf_{T\rightarrow\infty}\frac{M_{T}}{\sqrt{2\braket{M}_{T}\log\log\braket{M}_{T}}}=-1,\quad a.s.
				\end{aligned}
			\end{equation*}
		\end{enumerate}
	\end{prop}
	\begin{rmk}
		For the proof of (1) of Proposition \ref{propset}, one can refer to Theorem 1.6, Chapter 5, p.188 of \cite{Revuz and Yor(1998)}. While (2) of Proposition \ref{propset} {follows from Corollary 1, Chapter 2, p.144 of \cite{Liptser and Shiryayev(2011)}.}
	 And (3) of Proposition \ref{propset} can be obtained by the  distribution of Brownian motion, i.e., $B_{T}/\sqrt{T}\sim\mathcal{N}(0,1)$, for all $T>0$.
		Furthermore, (4) of Proposition \ref{propset} is obtained by the law of the iterated logarithm \citep[see Theorem 11.1, Section 11.1, p.165 of][]{Schilling et al.(2012)} of Brownian motion, i.e.,
		\begin{equation*}
			\begin{aligned}
				&\varlimsup_{T\rightarrow\infty}\frac{W_{T}}{\sqrt{2T\log\log T}}=1,\\ 
				&\varliminf_{T\rightarrow\infty}\frac{W_{T}}{\sqrt{2T\log\log T}}=-1, \quad a.s.
			\end{aligned}
		\end{equation*}
\end{rmk}

\section{The LSE based on continuous sampling}\label{secc}
In this section, we construct the continuous-type LSEs for the parameters $\alpha_{j}$ and $\beta_{j}$, $ j=1, \cdots, m$, and prove their consistency and asymptotic normality {based on the continuous observations.} We separate the discussion into two subsections: the first subsection considers that the drift term is known to be piecewise linear with $\beta_{j}=0$, while the second subsection considers that the drift term is known to be piecewise linear with unknown and {not vanishing piecewise constants $\beta_{j}$}.

\subsection{Estimate $\alpha_{j}$ for known $\beta_{j}=0$ and $\Theta$ $(j=1, \cdots, m)$}\label{secc1}
In this case, the SDE becomes 
\begin{equation}\label{eqmodelc1}
	\left\{
	\begin{aligned}
		&\mathrm{d}X_{t}=\sum_{j=1}^{m}-\alpha_{j}X_{t}I_{\{X_{t},j\}}\mathrm{d}t+\sigma \mathrm{d}W_{t},\\
		&X_{0}=x.
	\end{aligned}
	\right.
\end{equation}
When there is only one regime, i.e., $X$ is a standard OU process, the continuous-type LSE for $\alpha$ is given by minimizing the following objective function \citep[see for example][and
	the references therein]{Kutoyants(2004), Hu and Nualart(2010)}
	\begin{equation*}
		\int_{0}^{T}\left\lvert\dot{X}_{t}+\alpha X_{t}\right\rvert^{2}\mathrm{d}t.
	\end{equation*}
	The minimum is achieved at 
	\begin{equation*}
		\hat{\alpha}_{T}=-\frac{\int_{0}^{T}X_{t}\mathrm{d}X_{t}}{\int_{0}^{T}X_{t}^{2}\mathrm{d}t}.
	\end{equation*}
	When there are multiple regimes,
	the continuous-type LSE for $\alpha_{j}$ is motivated by the following heuristic argument, which aims to minimize
	\begin{equation*}
		\int_{0}^{T}\bigg|\dot{X}_{t}+\sum_{j=1}^{m}\alpha_{j}X_{t}I_{\{X_{t},j\}}\bigg|^{2}\mathrm{d}t,
	\end{equation*}
	where $\dot{X}_{t}$ denotes the formal derivative of $X_{t}$ with respect to $t$. 
	This is a quadratic function of $\alpha_{i}$, $i=1, \cdots, m$, although $\dot{X}_{t}$ does not exist. The minimum is achieved at 
	\begin{equation}\label{eqc1estimator}
		\hat{\alpha}_{i,T}=-\frac{\int_{0}^{T}X_{t}I_{\{X_{t},i\}}\left(\mathrm{d}X_{t}+\sum_{j\neq i}\alpha_{j}X_{t}I_{\{X_{t},j\}}\mathrm{d}t\right)}{\int_{0}^{T}X^{2}_{t}I_{\{X_{t},i\}}\mathrm{d}t}.
	\end{equation}
	By Eq. (\ref{eqmodelc1}), the facts that $I_{\{X_{t},i\}}I_{\{X_{t},i\}}=1$ and $I_{\{X_{t},i\}} I_{\{X_{t},j\}}=0$, for $i\neq j$, we have
	\begin{equation*}
		\hat{\alpha}_{i,T}=-\frac{\int_{0}^{T}X_{t}I_{\{X_{t},i\}}\mathrm{d}X_{t}}{\int_{0}^{T}X^{2}_{t}I_{\{X_{t},i\}}\mathrm{d}t}=\alpha_{i}-\sigma\frac{\int_{0}^{T}X_{t}I_{\{X_{t},i\}}\mathrm{d}W_{t}}{\int_{0}^{T}X^{2}_{t}I_{\{X_{t},i\}}\mathrm{d}t}.
\end{equation*}
Recall the unique invariant density of $\{X_{t}\}_{t\geq0}$, we define
\begin{equation}\label{eqc1Sigmai}
	M_{i}=\int_{\theta_{i-1}}^{\theta_{i}}k_{i}x^{2}\exp\left(\frac{-\alpha_{i}x^{2}}{\sigma^{2}}\right)\mathrm{d}x, \quad \text{for}\quad i=1, \cdots, m.
\end{equation}

The asymptotic behavior of our estimators $\hat{\alpha}_{i,T}$ is given as follows.

\hfill

\begin{thm}\label{thmc1}
	For $i=1, 2, \cdots, m$, the continuous-type LSEs $\hat{\alpha}_{i,T}$ defined in Eq. (\ref{eqc1estimator}) have the following statistical properties.
	\begin{enumerate}[label*=(\arabic*)]
		\item The continuous-type LSEs $\hat{\alpha}_{i,T}$ admit the strong consistency, i.e.,
		\begin{equation*}
			\hat{\alpha}_{i,T}\stackrel{a.s.}{\longrightarrow}\alpha_{i}, \quad\text{as}\quad T\rightarrow\infty.
		\end{equation*}
		\item The continuous-type LSEs $\hat{\alpha}_{i,T}$ admit the asymptotic normality, i.e.,
		\begin{equation*}
			\sqrt{T}(\hat{\alpha}_{i,T}-\alpha_{i})\sim \mathcal{N}(0, \sigma^{2}M_{i}^{-1}), \quad\text{as}\quad T\rightarrow\infty,
		\end{equation*}
		where $M_{i}$ is defined in Eq. (\ref{eqc1Sigmai}).
	\end{enumerate}
\end{thm}

\hfill

\noindent\textbf{Proof of Theorem \ref{thmc1}}.
Let
	\begin{equation*}
		\hat{\alpha}_{i,T}=\alpha_{i}-\sigma\frac{G_{i,T}}{F_{i,T}},
	\end{equation*}
	where 
	\begin{equation*}
		\begin{aligned}
				G_{i,t}&=\int_{0}^{t}X_{s}I_{\{X_{s},i\}}\mathrm{d}W_{s}\quad\text{and}\quad
		F_{i,t}=\int_{0}^{t}X^{2}_{s}I_{\{X_{s},i\}}\mathrm{d}s.
		\end{aligned}
	\end{equation*}
	To prove the strong consistency, it suffices to show that $G_{i,T}/F_{i,T}$ converges to zero almost surely as $T\rightarrow\infty$. 
	{It is obvious that the quadratic variation process of $G_{i,t}$ is exactly $F_{i,t}$. By Eq. (\ref{eqergodic}), we have
	\begin{equation}\label{eqc1thm1}
		\lim_{T\rightarrow\infty}\frac{F_{i,T}}{T}=\lim_{T\rightarrow\infty}\frac{1}{T}\int_{0}^{T}X^{2}_{t}I_{\{X_{t},i\}}\mathrm{d}t=M_{i},\quad a.s.
	\end{equation}
Note that $G_{i}$ is a continuous local martingale.
Then by (2) of Proposition \ref{propset} and Eq. (\ref{eqc1thm1}), we conclude 
\begin{equation*}\label{eqc1thm4}
	\lim_{T\rightarrow\infty}\frac{G_{i,T}}{F_{i,T}}=0,\quad a.s.,
\end{equation*}
which completes the proof of the strong consistency. 

By (3) of Proposition \ref{propset}, we have
\begin{equation}\label{eqc1thm5}
	\frac{G_{i,T}}{\sqrt{F_{i,T}}}\sim\mathcal{N}(0,1),\quad\text{as}\quad T\rightarrow\infty.
\end{equation}
Furthermore, combining Eq. (\ref{eqc1thm1}) and (\ref{eqc1thm5}), and Slutsky's theorem yields
	\begin{equation*}
		\sqrt{T}\left(\hat{\alpha}_{i,T}-\alpha_{i}\right)\sim\mathcal{N}\left(0, \sigma^{2}M_{i}^{-1}\right),\quad\text{as}\quad T\rightarrow\infty.
	\end{equation*}
	Thus, we complete the proof of the asymptotic normality.}
\hfill$\square$

\subsection{Estimate $\alpha_{j}$ and $\beta_{j}$ for known $\Theta$ $(j=1,\cdots,m)$}\label{secc2}
In this case, $\beta_{j}\neq0$ and the SDE becomes 
\begin{equation}\label{eqmodelc2}
	\left\{
	\begin{aligned}
		&\mathrm{d}X_{t}=\sum_{j=1}^{m}(\beta_{j}-\alpha_{j}X_{t})I_{\{X_{t},j\}}\mathrm{d}t+\sigma \mathrm{d}W_{t},\\
		&X_{0}=x.
	\end{aligned}
	\right.
\end{equation}
The continuous-type LSEs for $\alpha_{j}$ and $\beta_{j}$ aim to minimize the following objective function
\begin{equation*}
	\int_{0}^{T}\bigg|\dot{X}_{t}-\sum_{j=1}^{m}\left(\beta_{j}-\alpha_{j}X_{t}\right)I_{\{X_{t},j\}}\bigg|^{2}\mathrm{d}t,
\end{equation*}
which is a quadratic function of $\alpha_{i}$ and $\beta_{i}$, $i=1, \cdots, m$. By some direct calculations, the minimum is achieved at 
\begin{equation}\label{eqc21estimator}
\begin{aligned}
	\hat{\alpha}_{i,T}=&\frac{\int_{0}^{T}I_{\{X_{t},i\}}\mathrm{d}X_{t}\int_{0}^{T} X_{t}I_{\{X_{t},i\}} \mathrm{d} t-L_{i,T} \int_{0}^{T} X_{t}I_{\{X_{t},i\}} \mathrm{d} X_{t}}{L_{i,T} \int_{0}^{T} X_{t}^{2}I_{\{X_{t},i\}} \mathrm{d} t-\left(\int_{0}^{T} X_{t}I_{\{X_{t},i\}} \mathrm{d} t\right)^{2}}\\ 
\hat{\beta}_{i,T}=&\frac{\int_{0}^{T}I_{\{X_{t},i\}}\mathrm{d}X_{t} \int_{0}^{T} X_{t}^{2}I_{\{X_{t},i\}} \mathrm{d} t-\int_{0}^{T} X_{t}I_{\{X_{t},i\}} \mathrm{d} X_{t} \int_{0}^{T} X_{t}I_{\{X_{t},i\}} \mathrm{d} t}{L_{i,T} \int_{0}^{T} X_{t}^{2} I_{\{X_{t},i\}}\mathrm{d} t-\left(\int_{0}^{T} X_{t}I_{\{X_{t},i\}} \mathrm{d} t\right)^{2}},
\end{aligned}
\end{equation}
where $L_{i,T}=\int_{0}^{T}I_{\{X_{t},i\}}\mathrm{d}t$.
For simplicity, we make some notations as follows
\begin{equation}\label{eqc2thm0}
	\begin{aligned}
		P_{i}=&\int_{\theta_{i-1}}^{\theta_{i}}\pi(x)\mathrm{d}x=\mathbb{P}\left(\theta_{i-1}\leq X_{t}<\theta_{i}\right),\\
		K_{i}=&\int_{\theta_{i-1}}^{\theta_{i}}k_{i}x\exp\left(\frac{-\alpha_{i}x^{2}+2\beta_{i}x}{\sigma^{2}}\right)\mathrm{d}x,\\
		R_{i}=&\int_{\theta_{i-1}}^{\theta_{i}}k_{i}x^{2}\exp\left(\frac{-\alpha_{i}x^{2}+2\beta_{i}x}{\sigma^{2}}\right)\mathrm{d}x.
	\end{aligned}
\end{equation}

The asymptotic behavior of the proposed estimators $\hat{\alpha}_{i,T}$ and $\hat{\beta}_{i,T}$ is demonstrated in the following theorem.

\hfill

\begin{thm}\label{thmc2}
	For $i=1, \cdots, m$, the continuous-type LSEs $\hat{\alpha}_{i,T}$ and $\hat{\beta}_{i,T}$ defined in Eq. (\ref{eqc21estimator}) have the following statistical properties.
	\begin{enumerate}[label*=(\arabic*)]
		\item The continuous-type LSEs $\hat{\alpha}_{i,T}$ and $\hat{\beta}_{i,T}$ admit the strong consistency, i.e.,
		\begin{equation*}
			\left(\hat{\alpha}_{i,T}, \hat{\beta}_{i,T}\right)\stackrel{a.s.}{\longrightarrow}\left(\alpha_{i}, \beta_{i}\right), \quad\text{as}\quad T\rightarrow\infty.
		\end{equation*}
		\item The continuous-type LSEs $\hat{\alpha}_{i,T}$ and $\hat{\beta}_{i,T}$ admit the asymptotic normality, i.e.,
		\begin{equation*}
			\begin{aligned}
				&\sqrt{T}\left(\hat{\alpha}_{i,T}-\alpha_{i}\right)\sim \mathcal{N}\left(0,\frac{\sigma^{2}P_{i}}{ P_{i}R_{i}-K_{i}^{2}}\right),\\
				&\sqrt{T}\left(\hat{\beta}_{i,T}-\beta_{i}\right)\sim \mathcal{N}\left(0,\frac{\sigma^{2}R_{i}}{ P_{i}R_{i}-K_{i}^{2}}\right),
			\end{aligned}
		\end{equation*}
		as $T\rightarrow\infty$, where $P_{i}$, $R_{i}$ and $K_{i}$ are defined {in Eq. (\ref{eqc2thm0}). Furthermore}, the joint distribution of the continuous-type LSEs $\hat{\alpha}_{i,T}$ and $\hat{\beta}_{i,T}$ is given by
			{\begin{equation*}
				\sqrt{T}\left(\begin{array}{c}
					\hat{\alpha}_{i,T}-\alpha_{i}\\
					\hat{\beta}_{i,T}-\beta_{i}
				\end{array}\right)\sim\mathcal{N}\left[\left(\begin{array}{c}
					0\\
					0
				\end{array}\right), \left(
				\begin{array}{cc}
					\frac{\sigma^{2}P_{i}}{P_{i}R_{i}-K_{i}^{2}}, & \frac{\sigma^{2}K_{i}}{P_{i}R_{i}-K_{i}^{2}}\\
					\frac{\sigma^{2}K_{i}}{P_{i}R_{i}-K_{i}^{2}}, & \frac{\sigma^{2}R_{i}}{P_{i}R_{i}-K_{i}^{2}}
				\end{array}\right)
				\right],\quad\text{as}\quad T\rightarrow\infty.
		\end{equation*}}
	\end{enumerate}
\end{thm}

\hfill

\noindent\textbf{Proof of Theorem \ref{thmc2}}. By Eq. (\ref{eqmodelc2}), we have
\begin{equation}\label{eqc2thm1}
		\begin{aligned}
	&\int_{0}^{T}I_{\{X_{t},i\}}\mathrm{d}X_{t}=\beta_{i}L_{i,T}-\alpha_{i}\int_{0}^{T}X_{t}I_{\{X_{t},i\}}\mathrm{d}t+\sigma\int_{0}^{T}I_{\{X_{t},i\}}\mathrm{d}W_{t},\\
	&\int_{0}^{T}X_{t}I_{\{X_{t},i\}}\mathrm{d}X_{t}= \beta_{i}\int_{0}^{T}X_{t}I_{\{X_{t},i\}}\mathrm{d}t-\alpha_{i}\int_{0}^{T}X_{t}^{2}I_{\{X_{t},i\}}\mathrm{d}t+\sigma\int_{0}^{T}X_{t}I_{\{X_{t},i\}}\mathrm{d}W_{t}.
	\end{aligned}
\end{equation}
Hence, substituting Eq. (\ref{eqc21estimator})(\ref{eqc2thm1}) into Eq. (\ref{eqc2thm1}) yields
\begin{equation*}
	\begin{aligned}
		\hat{\alpha}_{i,T}&=\alpha_{i}+\sigma\frac{G_{i,T}}{F_{i,T}}\quad\text{and}\quad
		\hat{\beta}_{i,T}
		=\beta_{i}+\sigma\frac{H_{i,T}}{F_{i,T}},
	\end{aligned}
\end{equation*}
where 
\begin{equation*}
	\begin{aligned}
		G_{i,T}=&\frac{1}{T}\left(\int_{0}^{T}X_{t}I_{\{X_{t},i\}}\mathrm{d}t\int_{0}^{T}I_{\{X_{t},i\}}\mathrm{d}W_{t}-L_{i,T}\int_{0}^{T}X_{t}I_{\{X_{t},i\}}\mathrm{d}W_{t}\right),\\
		H_{i,T}=&\frac{1}{T}\left(\int_{0}^{T} X_{t}^{2}I_{\{X_{t},i\}} \mathrm{d} t\int_{0}^{T}I_{\{X_{t},i\}}\mathrm{d}W_{t} -\int_{0}^{T} X_{t}I_{\{X_{t},i\}} \mathrm{d} t\int_{0}^{T} X_{t}I_{\{X_{t},i\}} \mathrm{d} W_{t} \right),\\
		F_{i,T}=&\frac{1}{T}\left(L_{i,T} \int_{0}^{T} X_{t}^{2}I_{\{X_{t},i\}} \mathrm{d} t-\left(\int_{0}^{T} X_{t}I_{\{X_{t},i\}} \mathrm{d} t\right)^{2}\right).
	\end{aligned}
\end{equation*}
To prove the consistency of $\hat{\alpha}_{i,T}$ and $\hat{\beta}_{i,T}$, it suffices to prove that $G_{i,T}/F_{i,T}$ and $H_{i,T}/F_{i,T}$ converge to zero almost surely. By Eq. (\ref{eqergodic}), we have
\begin{equation}\label{eqc2thm3}
	\begin{aligned}
		&\lim_{T\rightarrow\infty}\frac{L_{i,T}}{T}=\int_{\theta_{i-1}}^{\theta_{i}}\pi(x)\mathrm{d}x=\mathbb{P}\left(\theta_{i-1}\leq X_{t}<\theta_{i}\right)=P_{i},\\
		&\lim_{T\rightarrow\infty}\frac{1}{T}\int_{0}^{T}X_{t}I_{\{X_{t},i\}}\mathrm{d}t=\int_{\theta_{i-1}}^{\theta_{i}}k_{i}x\exp\left(\frac{-\alpha_{i}x^{2}+2\beta_{i}x}{\sigma^{2}}\right)\mathrm{d}x=K_{i},\\
		&\lim_{T\rightarrow\infty}\frac{1}{T}\int_{0}^{T}X_{t}^{2}I_{\{X_{t},i\}}\mathrm{d}t=\int_{\theta_{i-1}}^{\theta_{i}}k_{i}x^{2}\exp\left(\frac{-\alpha_{i}x^{2}+2\beta_{i}x}{\sigma^{2}}\right)\mathrm{d}x=R_{i},\quad a.s.
	\end{aligned}
\end{equation}
On the one hand, by Eq. (\ref{eqc2thm3}), we have 
	\begin{equation}\label{eqc2thm4}
		\lim_{T\rightarrow\infty}\frac{F_{i,T}}{T}=P_{i}R_{i}-K_{i}^{2},\quad a.s.
\end{equation}
By It$\hat{\rm{o}}$'s isometry, the dominated convergence theorem, and Eq. (\ref{eqc2thm3}), we have
{
		\begin{align*}
			\lim_{T\rightarrow\infty}\frac{\braket{G_{i}}_{T}}{T}	=&K_{i}^{2}\lim_{T\rightarrow\infty}\frac{1}{T}\int_{0}^{T}I_{\{X_{t},i\}}\mathrm{d}t+P_{i}^{2}\lim_{T\rightarrow\infty}\frac{1}{T}\int_{0}^{T}X_{t}^{2}I_{\{X_{t},i\}}\mathrm{d}t\\&-2P_{i}K_{i}\lim_{T\rightarrow\infty}\frac{1}{T}\int_{0}^{T}X_{t}I_{\{X_{t},i\}}\mathrm{d}t\\
			=&P_{i}^{2}R_{i}-K_{i}^{2}P_{i},\quad a.s.
		\end{align*}
	and 
	\begin{equation*}
		\begin{aligned}
			\lim_{T\rightarrow\infty}\frac{\braket{H_{i}}_{T}}{T}
			=&R_{i}^{2}\lim_{T\rightarrow\infty}\frac{1}{T}\int_{0}^{T}I_{\{X_{t},i\}}\mathrm{d}t+K_{i}^{2}\lim_{T\rightarrow\infty}\frac{1}{T}\int_{0}^{T}X^{2}_{t}I_{\{X_{t},i\}}\mathrm{d}t\\&-2R_{i}K_{i}\lim_{T\rightarrow\infty}\frac{1}{T}\int_{0}^{T}X_{t}I_{\{X_{t},i\}}\mathrm{d}t\\
			=&R_{i}^{2}P_{i}-K_{i}^{2}R_{i},\quad a.s.
		\end{aligned}
	\end{equation*}}
 {Hence, we have
 \begin{equation*}
 \begin{aligned}
 &\lim_{T\rightarrow\infty}\braket{G_{i}}_{T}=\lim_{T\rightarrow\infty}\left(P_{i}^{2}R_{i}-K_{i}^{2}P_{i}\right)T=\infty,\\
 &\lim_{T\rightarrow\infty}\braket{H_{i}}_{T}=\lim_{T\rightarrow\infty}\left(R_{i}^{2}P_{i}-K_{i}^{2}R_{i}\right)T=\infty.
 \end{aligned}
 \end{equation*}}
	By (2) of Proposition \ref{propset}, we have
	\begin{equation}\label{eqc2thm61}
 \begin{aligned}
 &\lim_{T\rightarrow\infty}\frac{G_{i,T}}{\braket{G_{i}}_{T}}=0\quad\text{and}\quad
 \lim_{T\rightarrow\infty}\frac{H_{i,T}}{\braket{H_{i}}_{T}}=0,\quad a.s.
 \end{aligned}
	\end{equation}
	 Hence, we have
	\begin{equation}\label{eqc2thm62}
		\lim_{T\rightarrow\infty}\frac{G_{i,T}}{F_{i,T}}=0\quad\text{and}\quad
		\lim_{T\rightarrow\infty}\frac{H_{i,T}}{F_{i,T}}=0,\quad a.s.
	\end{equation}
 This completes the proof of strong consistency.
	
	Furthermore, by (3) of Proposition \ref{propset}, we have
	\begin{equation}\label{eqc2thm6}
		\frac{G_{i,T}}{\sqrt{\braket{G_{i}}_{T}}}\sim\mathcal{N}(0,1)\quad\text{and}\quad
		\frac{H_{i,T}}{\sqrt{\braket{H_{i}}_{T}}}\sim\mathcal{N}(0,1),
	\end{equation}
	as $T\rightarrow\infty$.
	Combining Eq. (\ref{eqc2thm4})-(\ref{eqc2thm6}), and Slutsky's theorem yields
	\begin{equation*}
	\begin{aligned}
	&\frac{\sigma\sqrt{T}G_{i,T}}{F_{i,T}}\sim\mathcal{N}\left(0,\frac{\sigma^{2}P_{i}}{ P_{i}R_{i}-K_{i}^{2}}\right),\\
	&\frac{\sigma\sqrt{T}H_{i,T}}{F_{i,T}}\sim\mathcal{N}\left(0,\frac{\sigma^{2}R_{i}}{ P_{i}R_{i}-K_{i}^{2}}\right),
	\end{aligned}
	\end{equation*}
	as $T\rightarrow\infty$, which implies
	\begin{equation}\label{eqalpha}
		\begin{aligned}
			&\sqrt{T}\left(\hat{\alpha}_{i,T}-\alpha_{i}\right)\sim \mathcal{N}\left(0,\frac{\sigma^{2}P_{i}}{ P_{i}R_{i}-K_{i}^{2}}\right),\\
			&\sqrt{T}\left(\hat{\beta}_{i,T}-\beta_{i}\right)\sim \mathcal{N}\left(0,\frac{\sigma^{2}R_{i}}{ P_{i}R_{i}-K_{i}^{2}}\right),
		\end{aligned}
	\end{equation}
	as $T\rightarrow\infty$.
	Thus, we complete the proof of the asymptotic normality.
	
	Now we consider the joint distribution of $\hat{\alpha}_{i,T}$ and $\hat{\beta}_{i,T}$. {By Eq. (\ref{eqc2thm4}) and the  dominated convergence theorem, the limit of the covariance for $\sqrt{T}(\hat{\alpha}_{i,T}-\alpha_{i})$ and $\sqrt{T}(\hat{\beta}_{i,T}-\beta_{i})$ is given by
	\begin{equation*}
		\begin{aligned}
			&\lim_{T\rightarrow\infty}Cov\left[\sqrt{T}(\hat{\alpha}_{i,T}-\alpha_{i}), \sqrt{T}(\hat{\beta}_{i,T}-\beta_{i})\right]\\
			=&\lim_{T\rightarrow\infty}\frac{\sigma^{2}}{\left(P_{i}R_{i}-K_{i}^{2}\right)^{2}T}\mathbb{E}\left[G_{i,T}H_{i,T}\right].
		\end{aligned}
	\end{equation*}
	By Eq. (\ref{eqc2thm3}) and the dominated convergence theorem again, we have
	\begin{equation*}
	\lim_{T\rightarrow\infty}\frac{\mathbb{E}\left[G_{i,T}H_{i,T}\right]}{T}=P_{i}K_{i}R_{i}-K_{i}^{3},\quad a.s.
	\end{equation*}}
	Hence, we have
	\begin{equation}\label{eqcov3}
		\begin{aligned}
			\lim_{T\rightarrow\infty}Cov\left[\sqrt{T}(\hat{\alpha}_{i,T}-\alpha_{i}), \sqrt{T}(\hat{\beta}_{i,T}-\beta_{i})\right]=\frac{\sigma^{2}K_{i}}{P_{i}R_{i}-K_{i}^{2}},\quad a.s.
		\end{aligned}
	\end{equation}
	Combining Eq. (\ref{eqalpha})-(\ref{eqcov3}) completes the proof.
	\hfill$\square$
	
	\hfill
	
	\begin{rmk}
		{It is obvious that $I_{\{X_{t},i\}}I_{\{X_{t},j\}}$ for $i\neq j$ and $t\geq0$. For all $T\geq0$, we have
		\begin{equation}\label{eqcov1}
  \begin{aligned}
  	&Cov\left[\sqrt{T}\left(\hat{\alpha}_{i,T}-\alpha_{i}\right),\sqrt{T}\left(\hat{\alpha}_{j,T}-\alpha_{j}\right)\right]=0,\\
   &Cov\left[\sqrt{T}\left(\hat{\alpha}_{i,T}-\alpha_{i}\right),\sqrt{T}\left(\hat{\beta}_{j,T}-\beta_{j}\right)\right]=0.
  \end{aligned}
		\end{equation}}
		From Eq. (\ref{eqcov1}) we have that $\sqrt{T}\left(\hat{\alpha}_{i,T}-\alpha_{i}\right)$ is independent of $\sqrt{T}\left(\hat{\alpha}_{j,T}-\alpha_{j}\right)$, and $\sqrt{T}\left(\hat{\alpha}_{i,T}-\alpha_{i}\right)$ is independent of $\sqrt{T}\left(\hat{\beta}_{j,T}-\beta_{j}\right)$ as $T\rightarrow\infty$.
\end{rmk}

\hfill

{\begin{rmk}
		The maximum likehood method and the least squares method yield the same objective function, as well as the formula of the estimators. Hence, the asymptotic behavior of the MLE and the LSE is also the same \citep{Su and Chan(2015),Mazzonetto and Pigato(2020)}. The two methods could deal with the cases of multiple thresholds with less calculations, which is a great improvement compared to the generalized moment methods.
\end{rmk}}

\section{The LSE based on discrete sampling}\label{secd}
 Assume that X is observed at regularly spaced time points $\{t_{k}:=kh\}_{k=0}^{n}$. In this section, we shall construct the discrete-type LSEs for the parameters $\alpha_{j}$ and $\beta_{j}$, $j=1,2,\cdots,m$ and prove their consistency and asymptotic normality based on the discretely observed process $\{X_{t_{k}}\}_{k=0}^{n}$. It is evident that one often encounters practical difficulties in obtaining a completely continuous observation of the sample path, and only the discretely observed process is possible. In this situation, one can approximate the (stochastic) integral by its ``Riemann-It$\hat{\rm{o}}$" sum to directly modify the continuous-type estimators to the discrete ones \citep[see, e.g.][]{Lejay and Pigato(2020),Le Breton(1976), Mazzonetto and Pigato(2020)}. However, we obtain the discrete-type estimators by constructing the contrast function. The discussion in this section is separated into two subsections as Section \ref{secc}.

To study the strong consistency and asymptotic normality of the LSEs for the drift parameters, we impose two hypotheses which propose the joint restrictions on the mesh size $h$ and the sample size $n$. Furthermore, $h=h_{n}$ depends on n. {

\hfill
\begin{hypo}\label{hypo1}
	$h\rightarrow0$ and $nh\rightarrow\infty$, as $n\rightarrow\infty$.
\end{hypo}

\hfill
\begin{hypo}\label{hypo2}
	$nh^{2}\log\log(1/h)\rightarrow0$, as $n\rightarrow\infty$.
\end{hypo}}

\hfill

The Hypothesis \ref{hypo1} guarantees that the strong consistency holds, while the Hypothesis \ref{hypo1}-\ref{hypo2} guarantee that the asymptotic normality holds.

\subsection{Estimate $\alpha_{j}$ for known $\beta_{j}=0$ and $\Theta$ $(j=1, \cdots, m)$}\label{secd1}
In this case, the threshold OU process is the solution to Eq. (\ref{eqmodelc1}).
To obtain the LSEs for $\alpha_{j}$, we introduce the following contrast function
\begin{equation*}
	\sum_{k=0}^{n-1}\bigg|X_{t_{k+1}}-X_{t_{k}}+\sum_{j=1}^{m}\alpha_{j}X_{t_{k}}I_{\{X_{t_{k}},j\}}h\bigg|^{2}.
\end{equation*}
This is a quadratic function of $\alpha_{i}$, $i=1,\cdots, m$, and the minimum is achieved at
\begin{equation}\label{eqd1estimator}
	\hat{\alpha}_{i,n}=-\frac{\sum_{k=0}^{n-1}(X_{t_{k+1}}-X_{t_{k}})X_{t_{k}}I_{\{X_{t_{k}},i\}}}{h\sum_{k=0}^{n-1}X^{2}_{t_{k}}I_{\{X_{t_{k}},i\}}},
\end{equation}
for $i=1, \cdots, m$.

The asymptotic behavior of our estimators $\hat{\alpha}_{i,n}$ is given as follows.

\hfill

\begin{thm}\label{thmd1}
	For $i=1, 2, \cdots, m$, the discrete-type LSEs $\hat{\alpha}_{i,n}$ defined in Eq. (\ref{eqd1estimator}) have the following statistical properties.
	\begin{enumerate}[label*=(\arabic*)]
		\item Under the Hypothesis \ref{hypo1}, the discrete-type LSEs $\hat{\alpha}_{i,n}$ admit the strong consistency, i.e.,
			\begin{equation*}\label{eqdis1}
				\hat{\alpha}_{i,n}\stackrel{a.s.}{\longrightarrow}\alpha_{i},\quad\text{as}\quad n\rightarrow\infty.
		\end{equation*}
		\item Under the Hypothesis \ref{hypo1}-\ref{hypo2}, the discrete-type LSEs $\hat{\alpha}_{i,n}$ admit the asymptotic normality, i.e.,
		\begin{equation*}
			\sqrt{nh}(\hat{\alpha}_{i,n}-\alpha_{i})\sim \mathcal{N}(0, \sigma^{2}M_{i}^{-1}), \quad\text{as}\quad n\rightarrow\infty,
		\end{equation*}
		where $M_{i}$ is defined in Eq. (\ref{eqc1Sigmai}).
	\end{enumerate}
\end{thm}

\hfill

\noindent\textbf{Proof of Theorem \ref{thmd1}}.
Note that
\begin{equation}\label{eqd1thm1}
	X_{t_{k+1}}-X_{t_{k}}=-\int_{t_{k}}^{t_{k+1}}\sum_{j=1}^{m}\alpha_{j}X_{t}I_{\{X_{t},j\}}\mathrm{d}t+\sigma \Delta W_{t_{k}},
\end{equation}
where $\Delta W_{t_{k}}=W_{t_{k+1}}-W_{t_{k}}$. Hence, we have
	\begin{align*}
		\hat{\alpha}_{i,n}=&-\frac{\sum_{k=0}^{n-1}(-\int_{t_{k}}^{t_{k+1}}\sum_{j=1}^{m}\alpha_{j}X_{t}I_{\{X_{t},j\}}\mathrm{d}t+\sigma\Delta W_{t_{k}})X_{t_{k}}I_{\{X_{t_{k}},i\}}}{h\sum_{k=0}^{n-1}X^{2}_{t_{k}}I_{\{X_{t_{k}},i\}}}\\
		=&\alpha_{i}+\frac{G_{i,n}}{F_{i,n}}-\sigma\frac{H_{i,n}}{F_{i,n}},
	\end{align*}
where
	\begin{align*}
		G_{i,n}=&\sum_{k=0}^{n-1}X_{t_{k}}I_{\{X_{t_{k}},i\}}\int_{t_{k}}^{t_{k+1}}\sum_{j=1}^{m}\alpha_{j}(X_{t}I_{\{X_{t},j\}}-X_{t_{k}}I_{\{X_{t_{k}},j\}})\mathrm{d}t,\\
		H_{i,n}=&\sum_{k=0}^{n-1} X_{t_{k}}I_{\{X_{t_{k}},i\}}\Delta W_{t_{k}},\quad	F_{i,n}=\sum_{k=0}^{n-1}X^{2}_{t_{k}}I_{\{X_{t_{k}},i\}}h.
	\end{align*}
Now, it suffices to show that $G_{i,n}/F_{i,n}$ and $H_{i,n}/F_{i,n}$ converge to zero almost surely. By Eq. (\ref{eqergodic}), we have
\begin{equation}\label{eqd1thm2}
	\lim_{n\rightarrow\infty}\frac{F_{i,n}}{nh}=M_{i},\quad a.s.,
\end{equation}
where $M_{i}$ is defined in Eq. (\ref{eqc1Sigmai}).
{Given $t\geq0$, $X_{t}$ is a Gaussian process with finite mean and finite variance. Additionally, the $p$-th moment of $X_{t}$ is finite for $p\geq0$, which implies $\sup_{t}|X_{t}|<\infty$, a.s.}
Let $\alpha^{\ast}=\max\{|\alpha_{1}|,\cdots,|\alpha_{m}|\}$.
By Eq. (\ref{eqd1thm1}), we have
	\begin{equation}\label{eqd1gron1}
		\begin{aligned}
			|X_{s}-X_{t_{k}}|\leq& \int_{t_{k}}^{s}\bigg|\sum_{j=1}^{m}\alpha_{j}\left(X_{t}I_{\{X_{t},j\}}-X_{t_{k}}I_{\{X_{t_{k}},j\}}\right)\bigg|\mathrm{d}t+\alpha^{\ast}|X_{t_{k}}|h+\sigma |\Delta W_{t_{k}}|\\
			\leq&\alpha^{\ast}\left(\int_{t_{k}}^{s}|X_{t}-X_{t_{k}}|\mathrm{d}t+\sup_{s\in(t_{k}, t_{k+1}]}|X_{s}|h+|X_{t_{k}}|h\right)+\sigma |\Delta W_{t_{k}}|\\
			\leq&\alpha^{\ast}\left(\int_{t_{k}}^{s}|X_{t}-X_{t_{k}}|\mathrm{d}t+2\sup_{s\in(t_{k}, t_{k+1}]}|X_{s}|h\right)+\sigma |\Delta W_{t_{k}}|.
		\end{aligned}
\end{equation}
{When $s=t_{k}$, we have that $|X_{s}-X_{t_{k}}|=0$.}
By Gronwall's inequality and the law of the iterated logarithm, we have
	{\begin{equation}\label{eqd1gron2}
		\begin{aligned}
			\sup_{s\in(t_{k}, t_{k+1}]}|X_{s}-X_{t_{k}}|\leq &\left(2\alpha^{\ast}\sup_{s\in(t_{k}, t_{k+1}]}|X_{s}|h+\sigma \sqrt{2h\log\log(1/h)}\right)e^{\alpha^{\ast}h}\\
			=&\mathcal{O}\left(\sqrt{h\log\log(1/h)}\right),\quad a.s.,
		\end{aligned}
\end{equation}
which goes to zero as $h\rightarrow0$. There exists a constant $C_{1}$ such that $X_{t}\in\left[X_{t_{k}}-C_{1}\sqrt{h\log\log(1/h)}, X_{t_{k}}+C_{1}\sqrt{h\log\log(1/h)}\right]$ for $t\in(t_{k},t_{k+1}]$. Without loss of generality, we could choose $C_{1}>\sqrt{2}\sigma$. Then we decompose $I_{\{X_{t_{k}},i\}}$ into three terms as follows
\begin{equation}\label{eqI123}
\begin{aligned}
	I_{\{X_{t_{k}},i\}}
	=I^{(1)}_{\{X_{t_{k}},i\}}+I^{(2)}_{\{X_{t_{k}},i\}}+I^{(3)}_{\{X_{t_{k}},i\}},
\end{aligned}
\end{equation}
where
\begin{equation*}
\begin{aligned}
	I^{(1)}_{\{X_{t_{k}},i\}}&=I_{\left\{\theta_{i-1}\leq X_{t_{k}}<\theta_{i-1}+C_{1}\sqrt{h\log\log(1/h)}\right\}},\\
	I^{(2)}_{\{X_{t_{k}},i\}}&=I_{\left\{\theta_{i-1}+C_{1}\sqrt{h\log\log(1/h)}\leq X_{t_{k}}< \theta_{i}-C_{1}\sqrt{h\log\log(1/h)}\right\}},\\
	I^{(3)}_{\{X_{t_{k}},i\}}&=I_{\left\{\theta_{i}-C_{1}\sqrt{h\log\log(1/h)}\leq X_{t_{k}}<\theta_{i}\right\}}.
\end{aligned}
\end{equation*}
For all $k\geq0$, we have
\begin{equation*}
\begin{aligned}
&\left\lvert X_{t_{k}}I^{(2)}_{\{X_{t_{k}},i\}}\int_{t_{k}}^{t_{k+1}}\sum_{j=1}^{m}\alpha_{j}(X_{t}I_{\{X_{t},j\}}-X_{t_{k}}I_{\{X_{t_{k}},j\}})\mathrm{d}t\right\rvert\\
\leq&\sup_{s\in(t_{k}, t_{k+1}]}|X_{s}|\int_{t_{k}}^{t_{k+1}}\alpha^{\ast}\sup_{s\in(t_{k}, t_{k+1}]}\lvert X_{s}-X_{t_{k}}\rvert\mathrm{d}t\\
\leq&\mathcal{O}\left(\sqrt{h^{3}\log\log(1/h)}\right),\quad a.s.
\end{aligned}
\end{equation*}
We thus have
\begin{equation*}
\begin{aligned}
&\lim_{n\rightarrow\infty}\frac{1}{nh}\left\lvert\sum_{k=0}^{n-1} X_{t_{k}}I^{(2)}_{\{X_{t_{k}},i\}}\int_{t_{k}}^{t_{k+1}}\sum_{j=1}^{m}\alpha_{j}(X_{t}I_{\{X_{t},j\}}-X_{t_{k}}I_{\{X_{t_{k}},j\}})\mathrm{d}t\right\rvert\\
=&\mathcal{O}\left(\sqrt{h\log\log(1/h)}\right),\quad a.s.
\end{aligned}
\end{equation*}
Furthermore, by Eq. (\ref{eqergodic}), we have
\begin{equation*}
\begin{aligned}
&\lim_{n\rightarrow\infty}\frac{1}{nh}\left\lvert\sum_{k=0}^{n-1} X_{t_{k}}\left(I^{(1)}_{\{X_{t_{k}},i\}}+I^{(3)}_{\{X_{t_{k}},i\}}\right)\int_{t_{k}}^{t_{k+1}}\sum_{j=1}^{m}\alpha_{j}(X_{t}I_{\{X_{t},j\}}-X_{t_{k}}I_{\{X_{t_{k}},j\}})\mathrm{d}t\right\rvert\\
\leq&\lim_{n\rightarrow\infty}\frac{1}{n}\sum_{k=0}^{n-1}\left(I^{(1)}_{\{X_{t_{k}},i\}}+I^{(3)}_{\{X_{t_{k}},i\}}\right)2\alpha^{\ast}\sup_{s\in(t_{k}, t_{k+1}]}\lvert X_{s}\rvert^{2}\\
\leq&2\alpha^{\ast}\sup_{t}|X_{t}|^{2}\left(\int_{\theta_{i-1}}^{\theta_{i-1}+C_{1}\sqrt{h\log\log(1/h)}}k_{i}e^{\frac{-\alpha_{i}x^{2}+2\beta_{i}x}{\sigma^{2}}}\mathrm{d}x\right.\\
&\left.+\int_{\theta_{i}-C_{1}\sqrt{h\log\log(1/h)}}^{\theta_{i}}k_{i}e^{\frac{-\alpha_{i}x^{2}+2\beta_{i}x}{\sigma^{2}}}\mathrm{d}x\right)\\
=&\mathcal{O}\left(\sqrt{h\log\log(1/h)}\right).
\end{aligned}
\end{equation*}
Hence, we have
\begin{equation}\label{eqd1thm3}
	\begin{aligned}
		\lim_{n\rightarrow\infty}\left\lvert\frac{G_{i,n}}{t_{n}}\right\rvert
		\leq &\mathcal{O}\left(\sqrt{h\log\log(1/h)}\right), \quad a.s.
	\end{aligned}
\end{equation}}
 Note that 
\begin{equation*}
	H_{i,m}=\sum_{k=0}^{m-1}\int_{t_{k}}^{t_{k+1}}X_{t_{k}}I_{\{X_{t_{k}},i\}}\mathrm{d}W_{t}\stackrel{a.s.}{\longrightarrow}\int_{0}^{t_{m}}X_{t}I_{\{X_{t},i\}}\mathrm{d}W_{t},\quad\text{as}\quad h\rightarrow0.
\end{equation*} 
Let $A_{i,t}=\int_{0}^{t}X_{t}I_{\{X_{t},i\}}\mathrm{d}W_{t}$ and $C_{i,t}$ be the quadratic variation process of $A_{i,t}$. {By Eq. (\ref{eqergodic}), we have
\begin{equation*}
	\begin{aligned}
		\lim_{n\rightarrow\infty}\frac{C_{i,t_{n}}}{nh}=\lim_{n\rightarrow\infty}\frac{1}{t_{n}}\int_{0}^{t_{n}}X_{t}^{2}I_{\{X_{t},i\}}\mathrm{d}t=M_{i}, \quad a.s.
	\end{aligned}
\end{equation*}}
{By (2) of Proposition \ref{propset}, we have}
\begin{equation}\label{eqdiscons2}
	\lim_{n\rightarrow\infty}\frac{H_{i,n}}{C_{i,t_{n}}}=\lim_{n\rightarrow\infty}\frac{H_{i,n}}{F_{i,n}}=0,\quad a.s.
\end{equation}
Combining Eq. (\ref{eqd1thm2}), (\ref{eqd1thm3}), and (\ref{eqdiscons2}) completes the proof of the strong consistency.


Note that
\begin{equation}\label{eqd1thm5}
	\lim_{n\rightarrow\infty}\left\lvert\frac{\sqrt{nh}G_{i,n}}{F_{i,n}}\right\rvert\leq \mathcal{O}\left(\sqrt{nh^{2}\log\log(1/h)}\right), \quad a.s.,
\end{equation} 
which goes to zero as $n\rightarrow\infty$. Now, it suffices to show that $\sqrt{nh}H_{i,n}/F_{i,n}$ converges in law to a centered normal distribution as $n\rightarrow\infty$. {By (3) of Proposition \ref{propset}, we have
\begin{equation*}
\frac{A_{i,t_{n}}}{\sqrt{C_{i,t_{n}}}}\sim\mathcal{N}(0,1),\quad \text{as}\quad n\rightarrow\infty,
\end{equation*}
which implies}
\begin{equation}\label{eqd1thm6}
	\frac{H_{i,n}}{\sqrt{F_{i,n}}}\sim\mathcal{N}(0,1),\quad\text{as}\quad n\rightarrow\infty.
\end{equation}
Combining Eq. (\ref{eqd1thm2}), (\ref{eqd1thm5}), (\ref{eqd1thm6}), and Slutsky's theorem completes the proof of the asymptotic normality.
\hfill$\square$

\hfill


\subsection{Estimate $\alpha_{j}$ and $\beta_{j}$ for known $\Theta$ $(j=1, \cdots, m)$}\label{secd2}
In this case, the threshold OU process is the solution to Eq. (\ref{eqmodelc2}).
To obtain the discrete-type LSEs for $\alpha_{j}$ and $\beta_{j}$, we introduce the following contrast function
\begin{equation*}
	\sum_{k=0}^{n-1}\bigg|X_{t_{k+1}}-X_{t_{k}}-\sum_{j=1}^{m}(\beta_{j}-\alpha_{j}X_{t_{k}})I_{\{X_{t_{k}},j\}}h\bigg|^{2}.
\end{equation*}
This is a quadratic function of $\alpha_{i}$ and $\beta_{i}$, $i=1, \cdots, m$, and the minimum is achieved at
\begin{equation*}
\begin{aligned}
	\hat{\alpha}_{i,n}=&\frac{Q_{i,n}J_{i,n}-L_{i,n}H_{i,n}}{L_{i,n}D_{i,n}h-J_{i,n}^{2}h},\\
\hat{\beta}_{i,n}=&\frac{Q_{i,n}D_{i,n}-H_{i,n}J_{i,n}}{L_{i,n}D_{i,n}h-J_{i,n}^{2}h},
\end{aligned}
\end{equation*}
where 
\begin{equation*}
	\begin{aligned}
		&L_{i,n}=\sum_{k=0}^{n-1}I_{\{X_{t_{k}},i\}},\quad Q_{i,n}=\sum_{k=0}^{n-1}(X_{t_{k+1}}-X_{t_{k}})I_{\{X_{t_{k}},i\}},\quad
		D_{i,n}=\sum_{k=0}^{n-1}X_{t_{k}}^{2}I_{\{X_{t_{k}},i\}},\\
		&J_{i,n}=\sum_{k=0}^{n-1}X_{t_{k}}I_{\{X_{t_{k}},i\}},\quad
		H_{i,n}=\sum_{k=0}^{n-1}X_{t_{k}}I_{\{X_{t_{k}},i\}}(X_{t_{k+1}}-X_{t_{k}}).
	\end{aligned}
\end{equation*}

The asymptotic behavior of the proposed estimators $\hat{\alpha}_{i,n}$ and $\hat{\beta}_{i,n}$ are given as follows.

\hfill

\begin{thm}\label{thmd2}
	For $i=1, \cdots, m$, the discrete-type LSEs $\hat{\alpha}_{i,n}$ and $\hat{\beta}_{i,n}$ defined in Eq. (\ref{eqc21estimator}) have the following statistical properties.
	\begin{enumerate}[label*=(\arabic*)]
		\item Under the Hypothesis \ref{hypo1}, the discrete-type LSEs $\hat{\alpha}_{i,n}$ and $\hat{\beta}_{i,n}$ admit the strong consistency, i.e.,
			\begin{equation*}\label{eqdis2}
				\left(\hat{\alpha}_{i,n}, \hat{\beta}_{i,n}\right)\stackrel{a.s.}{\longrightarrow}\left(\alpha_{i}, \beta_{i}\right),\quad\text{as}\quad n\rightarrow\infty.
		\end{equation*}
		\item Under the Hypothesis \ref{hypo1}-\ref{hypo2}, the discrete-type LSEs $\hat{\alpha}_{i,n}$ and $\hat{\beta}_{i,n}$ admit the asymptotic normality, i.e.,
		\begin{equation*}
			\begin{aligned}
				&\sqrt{nh}\left(\hat{\alpha}_{i,n}-\alpha_{i}\right)\sim \mathcal{N}\left(0, \frac{\sigma^{2}P_{i}}{ P_{i}R_{i}-K_{i}^{2}}\right),\\
				&\sqrt{nh}\left(\hat{\beta}_{i,n}-\beta_{i}\right)\sim \mathcal{N}\left(0,\frac{\sigma^{2}R_{i}}{ P_{i}R_{i}-K_{i}^{2}}\right),
			\end{aligned}
		\end{equation*}
		as $n\rightarrow\infty$, where $P_{i}$, $R_{i}$ and $K_{i}$ are defined in Eq. (\ref{eqc2thm3}). Furthermore, the joint distribution of the discrete-type LSEs $\hat{\alpha}_{i,n}$ and $\hat{\beta}_{i,n}$ is given by
			{\begin{equation*}
				\sqrt{nh}\left(\begin{array}{c}
					\hat{\alpha}_{i,n}-\alpha_{i}\\
					\hat{\beta}_{i,n}-\beta_{i}
				\end{array}\right)\sim\mathcal{N}\left[\left(\begin{array}{c}
					0\\
					0
				\end{array}\right), \left(
				\begin{array}{cc}
					\frac{\sigma^{2}P_{i}}{P_{i}R_{i}-K_{i}^{2}}, & \frac{\sigma^{2}K_{i}}{P_{i}R_{i}-K_{i}^{2}}\\
					\frac{\sigma^{2}K_{i}}{P_{i}R_{i}-K_{i}^{2}}, & \frac{\sigma^{2}R_{i}}{P_{i}R_{i}-K_{i}^{2}}
				\end{array}\right)
				\right],\quad\text{as}\quad n\rightarrow\infty.
		\end{equation*}}
	\end{enumerate}
\end{thm}

\hfill

\noindent\textbf{Proof of Theorem {\ref{thmd2}}}. Note that
\begin{equation}\label{eqDelta}
	X_{t_{k+1}}-X_{t_{k}}=\int_{t_{k}}^{t_{k+1}}\sum_{j=1}^{m}\left(\beta_{j}-\alpha_{j}X_{t}\right)I_{\{X_{t},j\}}\mathrm{d}t+\sigma\Delta W_{t_{k}}.
\end{equation}
Hence, we have
\begin{equation*}
	\begin{aligned}
		\hat{\alpha}_{i,n}
		=&\alpha_{i}+\frac{R_{i,n}+M_{i,n}+\sigma Z_{i,n}}{F_{i,n}},\\
		\hat{\beta}_{i,n}=&
		\beta_{i}+\frac{R^{\prime}_{i,n}+M^{\prime}_{i,n}+\sigma Z^{\prime}_{i,n}}{F_{i,n}},
	\end{aligned}
\end{equation*}
where
\begin{equation*}
	\begin{aligned}
		&R_{i,n}=-\frac{1}{n}L_{i,n}\sum_{k=0}^{n-1}X_{t_{k}}I_{\{X_{t_{k}},i\}}(p_{k}-q_{k}),\\
		&M_{i,n}=\frac{1}{n}J_{i,n}\sum_{k=0}^{n-1}I_{\{X_{t_{k}},i\}}(p_{k}-q_{k}),\\
		&Z_{i,n}=\frac{1}{n}J_{i,n}\sum_{k=0}^{n-1}I_{\{X_{t_{k}},i\}}\Delta W_{t_{k}}-\frac{1}{n}L_{i,n}\sum_{k=0}^{n-1}X_{t_{k}}I_{\{X_{t_{k}},i\}}\Delta W_{t_{k}},\\
		&F_{i,n}=\frac{1}{n}L_{i,n}D_{i,n}h-\frac{1}{n}J_{i,n}^{2}h,\\
		&R^{\prime}_{i,n}=\frac{1}{n}D_{i,n}\sum_{k=0}^{n-1}I_{\{X_{t_{k}},i\}}(p_{k}-q_{k}),\\
		&M^{\prime}_{i,n}=-\frac{1}{n}J_{i,n}\sum_{k=0}^{n-1}X_{t_{k}}I_{\{X_{t_{k}},i\}}(p_{k}-q_{k}),\\
		&Z^{\prime}_{i,n}=\frac{1}{n}D_{i,n}\sum_{k=0}^{n-1}I_{\{X_{t_{k}},i\}}\Delta W_{t_{k}}-\frac{1}{n}J_{i,n}\sum_{k=0}^{n-1}X_{t_{k}}I_{\{X_{t_{k}},i\}}\Delta W_{t_{k}},\\
		&p_{k}=\int_{t_{k}}^{t_{k+1}}\sum_{j=1}^{m}\beta_{j}\left(I_{\{X_{t},j\}}-I_{\{X_{t_{k}},j\}}\right)\mathrm{d}t,\\
		&q_{k}=\int_{t_{k}}^{t_{k+1}}\sum_{j=1}^{m}\alpha_{j}\left(X_{t}I_{\{X_{t},j\}}-X_{t_{k}}I_{\{X_{t_{k}},j\}}\right)\mathrm{d}t.
	\end{aligned}
\end{equation*}
{For all $t\geq0$, $X_{t}$ is a Gaussian random variable with finite mean and finite variance, which implies that $\sup_{t}|X_{t}|<\infty$, a.s.}
Let $\alpha^{\ast}=\max\{|\alpha_{1}|, \cdots, |\alpha_{m}|\}$ and $\beta^{\ast}=\max\{|\beta_{1}|, \cdots, |\beta_{m}|\}$. By some similar arguments as Eq. (\ref{eqd1gron1}) and (\ref{eqd1gron2}), we have 
	\begin{equation}\label{eqd2thm1}
		\begin{aligned}
			\sup_{s\in(t_{k}, t_{k+1}]}|X_{s}-X_{t_{k}}|\leq&\left(2\beta^{\ast}h+2\alpha^{\ast}\sup_{s\in(t_{k}, t_{k+1}]}|X_{s}|h+\sigma \sqrt{2h\log\log(1/h)}\right)e^{\alpha^{\ast}h}\\=& \mathcal{O}\left(\sqrt{h\log\log(1/h)}\right).
		\end{aligned}
\end{equation}
{It is non-trivial to show that $\frac{1}{nh}\sum_{k=0}^{n-1}I_{\{X_{t_{k}},i\}}q_{k}$, $\frac{1}{nh}\sum_{k=0}^{n-1}X_{t_{k}}I_{\{X_{t_{k}},i\}}q_{k}$, $\frac{1}{nh}\sum_{k=0}^{n-1}I_{\{X_{t_{k}},i\}}p_{k}$, and $\frac{1}{nh}\sum_{k=0}^{n-1}X_{t_{k}}I_{\{X_{t_{k}},i\}}p_{k}$ converge to zero. There exists a constant $C_{1}>\sqrt{2}\sigma$ such that $X_{t}\in\left[X_{t_{k}}-C_{1}\sqrt{h\log\log(1/h)}, X_{t_{k}}+C_{1}\sqrt{h\log\log(1/h)}\right]$ for $t\in(t_{k},t_{k+1}]$. By Eq. (\ref{eqd2thm1}), we decompose $I_{\{X_{t_{k}},i\}}$ into three terms as Eq. (\ref{eqI123}). Furthermore, we have
	\begin{equation*}
	\begin{aligned}
&I_{\{X_{t_{k}},i\}}q_{k}=\left(I^{(1)}_{\{X_{t_{k}},i\}}+I^{(2)}_{\{X_{t_{k}},i\}}+I^{(3)}_{\{X_{t_{k}},i\}}\right)q_{k},\\
&X_{t_{k}}I_{\{X_{t_{k}},i\}}q_{k}=X_{t_{k}}\left(I^{(1)}_{\{X_{t_{k}},i\}}+I^{(2)}_{\{X_{t_{k}},i\}}+I^{(3)}_{\{X_{t_{k}},i\}}\right)q_{k}.
	\end{aligned}
	\end{equation*}
For all $k$, we have
\begin{equation*}
\begin{aligned}
\left\lvert I^{(2)}_{\{X_{t_{k}},i\}}q_{k}\right\rvert=&\left\lvert I^{(2)}_{\{X_{t_{k}},i\}}\int_{t_{k}}^{t_{k+1}}\sum_{j=1}^{m}\alpha_{j}\left(X_{t}I_{\{X_{t},j\}}-X_{t_{k}}I_{\{X_{t_{k}},j\}}\right)\mathrm{d}t\right\rvert\\
\leq&\int_{t_{k}}^{t_{k+1}}\alpha_{i}\left\lvert X_{t}-X_{t_{k}}\right\rvert\mathrm{d}t\\
\leq&\alpha^{\ast}\sup_{t\in(t_{k},t_{k+1}]}|X_{t}-X_{t_{k}}|h\\
\leq&\mathcal{O}\left(\sqrt{h^{3}\log\log(1/h)}\right),\quad a.s.
\end{aligned}
\end{equation*}
which implies
\begin{equation*}
\lim_{n\rightarrow\infty}\frac{1}{nh}\left\lvert\sum_{k=0}^{n-1}I_{\{X_{t_{k}},i\}}^{(2)}q_{k}\right\rvert\leq\mathcal{O}\left(\sqrt{h\log\log(1/h)}\right), \quad a.s.
\end{equation*}
Furthermore, by Eq. (\ref{eqergodic}), we have
\begin{equation*}\label{eqqqq}
\begin{aligned}
&\lim_{n\rightarrow\infty}\frac{1}{nh}\left\lvert\sum_{k=0}^{n-1}\left(I_{\{X_{t_{k}},i\}}^{(1)}+I_{\{X_{t_{k}},i\}}^{(3)}\right)q_{k}\right\rvert\\
\leq&\lim_{n\rightarrow\infty}\frac{2\alpha^{\ast}}{n}\sum_{k=0}^{n-1}\left(I_{\{X_{t_{k}},i\}}^{(1)}+I_{\{X_{t_{k}},i\}}^{(3)}\right)\sup_{t\in(t_{k},t_{k+1}]}|X_{t}|\\
\leq&\lim_{n\rightarrow\infty}\frac{2\alpha^{\ast}\sup_{t}\lvert X_{t}\rvert}{n}\sum_{k=0}^{n-1}\left(I_{\{X_{t_{k}},i\}}^{(1)}+I_{\{X_{t_{k}},i\}}^{(3)}\right)\\
=&2\alpha^{\ast}\sup_{t\in[0,t_{n}]}\lvert X_{t}\rvert\left(\int_{\theta_{i-1}}^{\theta_{i-1}+C_{1}\sqrt{h\log\log(1/h)}}k_{i}e^{\frac{-\alpha_{i}x^{2}+2\beta_{i}x}{\sigma^{2}}}\mathrm{d}x\right.\\
&\left.+\int_{\theta_{i}-C_{1}\sqrt{h\log\log(1/h)}}^{\theta_{i}}k_{i}e^{\frac{-\alpha_{i}x^{2}+2\beta_{i}x}{\sigma^{2}}}\mathrm{d}x\right)\\
=&\mathcal{O}\left(\sqrt{h\log\log(1/h)}\right),\quad a.s.
\end{aligned}
\end{equation*}
From that $\sup_{t}|X_{t}|<\infty$, we have
\begin{align*}
\left\lvert X_{t_{k}}I_{\{X_{t_{k}},i\}}^{(2)}q_{k}\right\rvert\leq&\sup_{t}|X_{t}|\int_{t_{k}}^{t_{k+1}}\lvert\alpha_{i}\rvert|X_{t}-X_{t_{k}}|\mathrm{d}t\\
\leq&\alpha^{\ast}\sup_{t}|X_{t}|\sup_{t\in(t_{k}, t_{k+1}]}|X_{t}-X_{t_{k}}|h\\
\leq&\mathcal{O}\left(\sqrt{h^{3}\log\log(1/h)}\right),\quad a.s.
\end{align*}
which implies
\begin{equation*}
\lim_{n\rightarrow\infty}\frac{1}{nh}\left\lvert\sum_{k=0}^{n-1}X_{t_{k}}I_{\{X_{t_{k}},i\}}^{(2)}q_{k}\right\rvert\leq\mathcal{O}\left(\sqrt{h\log\log(1/h)}\right),\quad a.s.
\end{equation*}
Similarly, we have
\begin{equation*}
\begin{aligned}
&\lim_{n\rightarrow\infty}\frac{1}{nh}\left\lvert\sum_{k=0}^{n-1}X_{t_{k}}\left(I_{\{X_{t_{k}},i\}}^{(1)}+I_{\{X_{t_{k}},i\}}^{(3)}\right)q_{k}\right\rvert\\
\leq&\lim_{n\rightarrow\infty}\frac{1}{n}\sum_{k=0}^{n-1}\left(I_{\{X_{t_{k}},i\}}^{(1)}+I_{\{X_{t_{k}},i\}}^{(3)}\right)2\alpha^{\ast}\sup_{t\in(t_{k},t_{k+1}]}X_{t}^{2}\\
\leq&\lim_{n\rightarrow\infty}\frac{2\alpha^{\ast}\sup_{t}X_{t}^{2}}{n}\sum_{k=0}^{n-1}\left(I_{\{X_{t_{k}},i\}}^{(1)}+I_{\{X_{t_{k}},i\}}^{(3)}\right)\\
=&\mathcal{O}\left(\sqrt{h\log\log(1/h)}\right),\quad a.s.
\end{aligned}
\end{equation*}
Hence, we have
\begin{equation*}\label{eqd2thm2}
	\begin{aligned}
		&\lim_{n\rightarrow\infty}\frac{1}{nh}\left\lvert\sum_{k=0}^{n-1}I_{\{X_{t_{k}},i\}}q_{k}\right\rvert\leq \mathcal{O}\left(\sqrt{h\log\log(1/h)}\right),\quad a.s.\\
		&\lim_{n\rightarrow\infty}\frac{1}{nh}\left\lvert\sum_{k=0}^{n-1}X_{t_{k}}I_{\{X_{t_{k}},i\}}q_{k}\right\rvert\leq \mathcal{O}\left(\sqrt{h\log\log(1/h)}\right),\quad a.s.
	\end{aligned}
\end{equation*}
}
{For $s\in(t_{k},t_{k+1}]$ and $I^{(2)}_{\{X_{t_{k}},i\}}=1$, by Eq. (\ref{eqd2thm1}), we deduce
\begin{equation*}
	\sum_{j=1}^{m}\beta_{j}I_{\{X_{s},j\}}I^{(2)}_{\{X_{t_{k}},i\}}=\beta_{i}.
\end{equation*}}
Hence, we have
{\begin{equation*}
	\begin{aligned}
		I_{\{X_{t_{k}},i\}}p_{k}=\left(I^{(1)}_{\{X_{t_{k}},i\}}+I^{(3)}_{\{X_{t_{k}},i\}}\right)p_{k}\leq2\left(I^{(1)}_{\{X_{t_{k}},i\}}+I^{(3)}_{\{X_{t_{k}},i\}}\right)\beta^{\ast}h,\quad a.s.
	\end{aligned}
\end{equation*}}
Then, by Eq. (\ref{eqergodic}), we have
\begin{equation*}\label{eqd2thm3}
	\begin{aligned}
		\lim_{n\rightarrow\infty}\frac{1}{nh}\left|\sum_{k=0}^{n-1}X_{t_{k}}I_{\{X_{t_{k}},i\}}p_{k}\right|
		\leq&\lim_{n\rightarrow\infty}\frac{{2}}{n}\left|\sum_{k=0}^{n-1}X_{t_{k}}\left(I^{(1)}_{\{X_{t_{k}},i\}}+I^{(3)}_{\{X_{t_{k}},i\}}\right)\right|\beta^{\ast}\\
		=&\left|\int_{\theta_{i-1}}^{\theta_{i-1}+C_{1}\sqrt{h\log\log(1/h)}}k_{i}xe^{\frac{-\alpha_{i}x^{2}+2\beta_{i}x}{\sigma^{2}}}\mathrm{d}x\right.\\&\left.+\int_{\theta_{i}-C_{1}\sqrt{h\log\log(1/h)}}^{\theta_{i}}k_{i}xe^{\frac{-\alpha_{i}x^{2}+2\beta_{i}x}{\sigma^{2}}}\mathrm{d}x\right|\beta^{\ast}\\
		=&\mathcal{O}(\sqrt{h\log\log(1/h)}),\quad a.s.
	\end{aligned}
\end{equation*}
and
\begin{equation*}\label{eqd2thm4}
	\begin{aligned}
		\lim_{n\rightarrow\infty}\frac{1}{nh}\left|\sum_{k=0}^{n-1}I_{\{X_{t_{k}},i\}}p_{k}\right|
		\leq&{2}\left|\int_{\theta_{i-1}}^{\theta_{i-1}+C_{1}\sqrt{h\log\log(1/h)}}k_{i}e^{\frac{-\alpha_{i}x^{2}+2\beta_{i}x}{\sigma^{2}}}\mathrm{d}x\right.\\
		&\left.+\int_{\theta_{i}-C_{1}\sqrt{h\log\log(1/h)}}^{\theta_{i}}k_{i}e^{\frac{-\alpha_{i}x^{2}+2\beta_{i}x}{\sigma^{2}}}\mathrm{d}x\right|\beta^{\ast}\\
		=&\mathcal{O}(\sqrt{h\log\log(1/h)}),\quad a.s.
	\end{aligned}
\end{equation*}
Using Eq. (\ref{eqergodic}) again yields
\begin{equation*}\label{eqd2thm5}
	\lim_{n\rightarrow\infty}\frac{1}{n}L_{i,n}=P_{i}, \quad \lim_{n\rightarrow\infty}\frac{1}{n}D_{i,n}=R_{i},\quad \lim_{n\rightarrow\infty}\frac{1}{n}J_{i,n}=K_{i},\quad a.s.,
\end{equation*}
where $P_{i}$, $R_{i}$ and $K_{i}$ are defined in Eq. (\ref{eqc2thm3}). Hence, we have
\begin{equation}\label{eqd2thmF}
\lim_{n\rightarrow\infty}\frac{F_{i,n}}{nh}=P_{i}R_{i}-K_{i}^{2},\quad a.s.
\end{equation}
Furthermore, we have
\begin{equation}\label{eqd2R+M}
\begin{aligned}
	&\lim_{n\rightarrow\infty}\left\lvert\frac{R_{i,n}+M_{i,n}}{F_{i,n}}\right\rvert\leq \mathcal{O}\left(\sqrt{h\log\log(1/h)}\right), \quad a.s.,\\
&\lim_{n\rightarrow\infty}\left\lvert\frac{R^{\prime}_{i,n}+M^{\prime}_{i,n}}{F_{i,n}}\right\rvert\leq \mathcal{O}\left(\sqrt{h\log\log(1/h)}\right),\quad a.s.,
\end{aligned}
\end{equation}
which tend to $0$ as $h\rightarrow0$.
{Let $A_{i,t}=K_{i}\int_{0}^{t}I_{\{X_{s},i\}}\mathrm{d}W_{s}-P_{i}\int_{0}^{t}X_{s}I_{\{X_{s},i\}}\mathrm{d}W_{s}$ and $A^{\prime}_{i,t}=R_{i}\int_{0}^{t}I_{\{X_{s},i\}}\mathrm{d}W_{s}-K_{i}\int_{0}^{t}X_{s}I_{\{X_{s},i\}}\mathrm{d}W_{s}$. Then we have
\begin{equation*}
Z_{i,n}\stackrel{a.s.}{\longrightarrow}A_{i,t_{n}},\quad Z_{i,n}^{\prime}\stackrel{a.s.}{\longrightarrow}A_{i,t_{n}}^{\prime},\quad \text{as}\quad n\rightarrow\infty.
\end{equation*}
Let $C_{i,t}$ and $C^{\prime}_{i,t}$ be the quadratic variation process of $A_{i,t}$ and $A_{i,t}^{\prime}$ respectively. By Eq. (\ref{eqergodic}), we have
\begin{equation*}
\lim_{n\rightarrow\infty}\frac{C_{i,t_{n}}}{nh}=P_{i}^{2}R_{i}-K_{i}^{2}P_{i}\quad\text{and}\quad \lim_{n\rightarrow\infty}\frac{C^{\prime}_{i,t_{n}}}{nh}=R_{i}^{2}P_{i}-K_{i}^{2}R_{i}, \quad a.s.
\end{equation*}
By (2) of Proposition \ref{propset}, we have
\begin{equation*}
\lim_{n\rightarrow\infty}\frac{A_{i,t_{n}}}{C_{i,t_{n}}}=0\quad\text{and}\quad\lim_{n\rightarrow\infty}\frac{A^{\prime}_{i,t_{n}}}{C^{\prime}_{i,t_{n}}}=0,\quad a.s.
\end{equation*}
Furthermore, we have
\begin{equation}\label{eqd2thm7}
\begin{aligned}
&\lim_{n\rightarrow\infty}\frac{Z_{i,n}}{F_{i,n}}=\lim_{n\rightarrow\infty}\frac{A_{i,n}}{C_{i,n}}\frac{C_{i,n}}{F_{i,n}}=0,\quad a.s.\\
&\lim_{n\rightarrow\infty}\frac{Z^{\prime}_{i,n}}{F_{i,n}}=\lim_{n\rightarrow\infty}\frac{A^{\prime}_{i,n}}{C^{\prime}_{i,n}}\frac{C^{\prime}_{i,n}}{F_{i,n}}=0,\quad a.s.
\end{aligned}
\end{equation}}
By Eq. (\ref{eqd2R+M})-(\ref{eqd2thm7}), we complete the proof of the strong consistency.

%
Now, we are in a position to prove the asymptotic normality. By Eq. (\ref{eqd2R+M}), we have
\begin{equation}\label{eqd2thm0}
\begin{aligned}
&\lim_{n\rightarrow\infty}\left\lvert\frac{\sqrt{nh}\left(R_{i,n}+M_{i,n}\right)}{F_{i,n}}\right\rvert\leq \mathcal{O}\left(\sqrt{nh^{2}\log\log(1/h)}\right), \quad a.s.,\\
&\lim_{n\rightarrow\infty}\left\lvert\frac{\sqrt{nh}\left(R^{\prime}_{i,n}+M^{\prime}_{i,n}\right)}{F_{i,n}}\right\rvert\leq \mathcal{O}\left(\sqrt{nh^{2}\log\log(1/h)}\right),\quad a.s.,
\end{aligned}
\end{equation}
which tend to zero as $n\rightarrow\infty$.
{By (3) of Proposition \ref{propset}, we have
\begin{equation}\label{eqd2thmZ}
	\frac{Z_{i,n}}{\sqrt{\left(P_{i}^{2}R_{i}-K_{i}^{2}P_{i}\right)nh}}\sim\mathcal{N}(0,1)
	\quad\text{and}\quad 
	\frac{Z^{\prime}_{i,n}}{\sqrt{\left(R_{i}^{2}P_{i}-K_{i}^{2}R_{i}\right)nh}}\sim\mathcal{N}(0,1),
\end{equation}}
as $n\rightarrow\infty$. Combining Eq. (\ref{eqd2thmF}), (\ref{eqd2thm0})-(\ref{eqd2thmZ}), and Slutsky's theorem yields
	\begin{equation}\label{eqalphad}
		\begin{aligned}
			&\sqrt{nh}\left(\hat{\alpha}_{i,n}-\alpha_{i}\right)\sim \mathcal{N}\left(0, \frac{\sigma^{2}P_{i}}{P_{i}R_{i}-K_{i}^{2}}\right),\\
			&\sqrt{nh}\left(\hat{\beta}_{i,n}-\beta_{i}\right)\sim \mathcal{N}\left(0,\frac{\sigma^{2}R_{i}}{ P_{i}R_{i}-K_{i}^{2}}\right),
		\end{aligned}
	\end{equation}
	as $n\rightarrow\infty$.
	
	Now we are in a position to prove the joint distribution of $\hat{\alpha}_{i,n}$ and $\hat{\beta}_{i,n}$. {Similar arguments as those	to derive Eq. (\ref{eqcov3}) yield}
	\begin{equation}\label{eqdcov3}
		\begin{aligned}
			\lim_{n\rightarrow\infty}Cov\left[\sqrt{nh}(\hat{\alpha}_{i,n}-\alpha_{i}), \sqrt{nh}(\hat{\beta}_{i,n}-\beta_{i})\right]=\frac{\sigma^{2}K_{i}}{P_{i}R_{i}-K_{i}^{2}},\quad a.s.
		\end{aligned}
	\end{equation}
	Combining Eq. (\ref{eqalphad})-(\ref{eqdcov3}) completes the proof.
\hfill$\square$

\hfill
 
	\begin{rmk}
		{Let 
	\begin{equation*}
		Z^{(i)}_{t}=K_{i}\int_{0}^{t}I_{\{X_{s},i\}}\mathrm{d}W_{s}-P_{i}\int_{0}^{t}X_{s}I_{\{X_{s},i\}}\mathrm{d}W_{s}.
	\end{equation*}
  Then we have $Z_{i,n}\stackrel{a.s.}{\longrightarrow}Z^{(i)}_{t_{n}}$, as $h\rightarrow0$.
		By (4) of Proposition \ref{propset}, we have 
		\begin{equation*}
		\begin{aligned}
		&\varlimsup_{n\rightarrow\infty}\frac{Z_{i,n}}{\sqrt{2\braket{Z^{(i)}}_{t_{n}}\log\log\braket{Z^{(i)}}_{t_{n}}}}=1,\\
		&\varliminf_{n\rightarrow\infty}\frac{Z_{i,n}}{\sqrt{2\braket{Z^{(i)}}_{t_{n}}\log\log\braket{Z^{(i)}}_{t_{n}}}}=-1,\quad a.s.
		\end{aligned}
		\end{equation*}}
		Hence, we have
		\begin{equation}\label{eqorder}
			\lim_{n\rightarrow\infty}\left\lvert\frac{Z_{i,n}}{F_{i,n}}\right\rvert\leq \mathcal{O}\left(\sqrt{\frac{\log\log nh}{nh}}\right),\quad a.s.
		\end{equation}
		Combining Eq. (\ref{eqd2R+M}) and (\ref{eqorder}) yields
		\begin{equation}\label{eqalphao}
			\lim_{n\rightarrow\infty}\lvert\hat{\alpha}_{i,n}-\alpha_{i}\rvert\leq \mathcal{O}\left(\sqrt{\frac{\log\log nh}{nh}}\right),\quad a.s.
		\end{equation}
		Similarly, we have
		\begin{equation}\label{eqbetao}
			\lim_{n\rightarrow\infty}\lvert\hat{\beta}_{i,n}-\beta_{i}\rvert\leq \mathcal{O}\left(\sqrt{\frac{\log\log nh}{nh}}\right),\quad a.s.
		\end{equation}
	\end{rmk}

\hfill

{\begin{rmk}
When the observations are discrete, \cite{Mazzonetto and Pigato(2020)} propose the discretized likehood function by approximating the (stochastic) integral by its ``Riemann-It$\hat{\rm{o}}$'' sum. The least squares estimator is proprosed by discretizing the SDE and minimizing the sums of squares of error. The two methods yield the same formula for estimators, as well as the asymptotic behavior.
\end{rmk}}

\section{Parameter estimation for $\sigma$}\label{secsigma}
In \cite{Lejay and Pigato(2018)}, they study the asymptotic behavior of estimators for the oscillating Brownian motion, i.e., a two-valued, discontinuous diffusion coefficient of a SDE as follows
\begin{equation}\label{eqsigma2}
	\left\{
	\begin{aligned}
		&\mathrm{d}X_{t}=\sigma^{+}I_{\{X_{t}\geq 0\}}\mathrm{d}W_{t}+\sigma^{-}I_{\{X_{t}<0\}}\mathrm{d}W_{t},\\
		&X_{0}=x,
	\end{aligned}
	\right.
\end{equation}
where $\sigma^{+}, \sigma^{-}>0$. 

Now, we consider a generalized SDE { of the form}
\begin{equation}\label{eqsigma1}
	\left\{
	\begin{aligned}
		&\mathrm{d}X_{t}=\sum_{j=1}^{m}\left\{\left(\beta_{j}-\alpha_{j}X_{t}\right)\mathrm{d}t+\sigma_{j}\mathrm{d}W_{t}\right\}I_{\{X_{t},j\}},\\
		&X_{0}=x,
	\end{aligned}
	\right.
\end{equation}
where $\sigma_{j}>0$, $j=1,\cdots,m$ are the diffusion parameters and the others keep the same as Eq. (\ref{eq1}). The LSEs for $\alpha_{j}$ and $\beta_{j}$ proposed in Section \ref{secc} and \ref{secd} are suitable in this case when long-time observations are available. Hence, in this section, we focus on the problem of proposing estimators for $\sigma_{j}$. The SDE (\ref{eqsigma2}) is a special case of the above SDE (\ref{eqsigma1}), i.e., the drift term is zero, $m=2$, and the only threshold is zero. Motivated by \cite{Lejay and Pigato(2018)}, we shall propose the QVEs for the diffusion parameters $\sigma_{i}$.

Recall that QVE stands for quadratic variation estimator and MQVE stands for modified quadratic variation estimator. Based on the continuously observed process $\{X_{t}\}_{0\leq t\leq\Lambda}$, we define the QVEs for $\sigma_{i}$ {as follows}
\begin{equation}\label{eqsigmaconti}
	\tilde{\sigma}_{i,\Lambda}=\sqrt{\frac{\braket{XI_{\{X,i\}}}_{\Lambda}}{L_{i,\Lambda}}},
\end{equation}
where $L_{i,\Lambda}=\int_{0}^{\Lambda}I_{\{X_{t},i\}}\mathrm{d}t$.  We emphasize that the QVEs only need short time observations, i.e., $0<\Lambda<\infty$.


Denote by $\ell_{Y,t}^{x}$ the local time at a point $x$, which represents the time spent by $Y$ at $x$ until $t$, which means
	{\begin{equation*}
		\ell_{Y,t}^{x}=\lim_{\varepsilon\rightarrow0}\frac{1}{2\varepsilon}\int_{0}^{t}I_{\{x-\varepsilon\leq Y_{s}\leq x+\varepsilon\}}\mathrm{d}\braket{Y}_{s}.
\end{equation*}}

The following theorem shows the strong consistency of the continuous type QVEs $\tilde{\sigma}_{i,\Lambda}$.

\hfill

\begin{thm}\label{thmsigma1}
	For all $0<\Lambda<\infty$, the continuous-type QVEs of $\sigma_{i}$ defined in Eq. (\ref{eqsigmaconti}) admit the strong consistency, i.e.,
	{
		$\tilde{\sigma}_{i,\Lambda}=\sigma_{i}$, a.s.
	}
\end{thm}

\hfill

\noindent\textbf{Proof of Theorem \ref{thmsigma1}}. Decompose the $m$ regimes into three cases where $\theta_{j-1}<\theta_{j}<X_{0}$, $\theta_{j-1}\leq X_{0}\leq\theta_{j}$, and $X_{0}<\theta_{j-1}<\theta_{j}$. We give the proof of the third case $X_{0}<\theta_{j-1}<\theta_{j}$ and the proof of the other two cases is similar.

Consider the following process
\begin{equation*}
	X^{(i)}_{t}=\left\{
	\begin{aligned}
		0,\quad\quad\quad&\text{if}\quad X_{t}\in(-\infty,\theta_{i-1}),\\
		X_{t}-\theta_{i-1},\quad&\text{if}\quad X_{t}\in[\theta_{i-1},\theta_{i}),\\
		\theta_{i}-\theta_{i-1},\quad&\text{if}\quad X_{t}\in[\theta_{i},\infty).
	\end{aligned}
	\right.
\end{equation*}
Hence, we have
\begin{equation*}
	\mathrm{d}X^{(i)}_{t}=\left\{
	\begin{aligned}
		0,\quad\quad&\text{if}\quad X_{t}\in(-\infty,\theta_{i-1}),\\
		\mathrm{d}X_{t},\quad&\text{if}\quad X_{t}\in[\theta_{i-1},\theta_{i}),\\
		0,\quad\quad&\text{if}\quad X_{t}\in[\theta_{i},\infty).
	\end{aligned}
	\right.
\end{equation*}
Using  Meyer-Tanaka's formula for $X_{t}^{(i)}$ yields
\begin{equation*}
	X_{t}^{(i)}=\int_{0}^{t}\mathrm{d}X_{t}^{(i)}+\ell^{0}_{X^{(i)},t}.
\end{equation*}
By the fact that $\mathrm{d}X_{t}^{(i)}=I_{\{X_{t},i\}}\mathrm{d}X_{t}$, we have
\begin{equation*}
	\begin{aligned}
		X_{t}^{(i)}=&\int_{0}^{t}I_{\{X_{t},i\}}\mathrm{d}X_{t}+\ell^{\theta_{i-1}}_{X,t}/2
		=\xi_{1,t}^{(i)}+\xi_{2,t}^{(i)}+\ell^{\theta_{i-1}}_{X,t}/2,
	\end{aligned}
\end{equation*}
where $\xi_{1,t}^{(i)}=\int_{0}^{t}(\beta_{i}-\alpha_{i}X_{t})I_{\{X_{t},i\}}\mathrm{d}t$ and $\xi_{2,t}^{(i)}=\sigma_{i}\int_{0}^{t}I_{\{X_{t},i\}}\mathrm{d}W_{t}$. 
Furthermore, we have
\begin{equation*}
	\begin{aligned}
		\braket{X^{(i)}}_{\Lambda}=&\braket{\xi_{1}^{(i)}}_{\Lambda}+\braket{\xi^{(i)}_{2}}_{\Lambda}+\braket{\ell^{\theta_{i-1}}_{X}}_{\Lambda}/4+2\braket{\xi^{(i)}_{1},\xi^{(i)}_{2}}_{\Lambda}+\braket{\xi^{(i)}_{1},\ell^{\theta_{i-1}}_{X}}_{\Lambda}+\braket{\xi^{(i)}_{2},\ell^{\theta_{i-1}}_{X}}_{\Lambda}\\
		=&-\braket{\xi_{1}^{(i)}}_{\Lambda}+\braket{\xi_{2}^{(i)}}_{\Lambda}-\braket{\ell^{\theta_{i-1}}_{X}}_{\Lambda}/4-\braket{\xi_{1}^{(i)},\ell^{\theta_{i-1}}_{X}}_{\Lambda}\\&+2\braket{\xi_{1}^{(i)},X^{(i)}}_{\Lambda}+\braket{\ell^{\theta_{i-1}}_{X},X^{(i)}}_{\Lambda}.
	\end{aligned}
\end{equation*}
Note that $\xi_{2,\Lambda}^{(i)}$ is a martingale with quadratic variation 
\begin{equation*}
	\braket{\xi_{2}^{(i)}}_{\Lambda}=\sigma_{i}^{2}\int_{0}^{\Lambda}I_{\{X_{t},i\}}\mathrm{d}t,
\end{equation*}
which implies that $\braket{\xi_{2}^{(i)}}_{\Lambda}/L_{i,\Lambda}=\sigma_{i}^{2}$.
Furthermore, the local time $\ell_{X}^{\theta_{i-1}}$, $\xi_{1}^{(i)}$, and $X^{(i)}$ are continuous, and $\xi_{1}^{(i)}$ is of finite variation. Thus, $\braket{\xi_{1}^{(i)}}_{\Lambda}$, $\braket{\xi_{1}^{(i)},\ell_{X}^{\theta_{i-1}}}_{\Lambda}$, $\braket{\xi_{1}^{(i)},X^{(i)}}_{\Lambda}$, $\braket{\ell^{\theta_{i-1}}_{X}}_{\Lambda}$, and $\braket{\ell_{X}^{\theta_{i-1}},X^{(i)}}_{\Lambda}$ equal to zero almost surely, which implies
\begin{equation*}
	\frac{\braket{XI_{\{X,i\}}}}{L_{i,\Lambda}}=\sigma_{i}^{2}.
\end{equation*}
Hence, we complete the proof.
\hfill$\square$

\hfill

For any finite time interval $[0,\Lambda]$, we consider the discretely observed process $\{X_{t_{k}}\}_{k=0}^{n}$ in $[0,\Lambda]$ given by $t_{k}=kh$ and $h=\Lambda/n$. Let $[X,Y]_{\Lambda}^{n}=\sum_{k=0}^{n-1}(X_{t_{k+1}}-X_{t_{k}})(Y_{t_{k+1}}-Y_{t_{k}})$ and $[X]_{\Lambda}^{n}=[X,X]_{\Lambda}^{n}$. Then we define the discrete-type QVE as follows
\begin{equation*}
	\tilde{\sigma}_{i,n}=\sqrt{\frac{[XI_{\{X,i\}}]_{\Lambda}^{n}}{L_{i,n}}},
\end{equation*}
where $L_{i,n}=\sum_{k=0}^{n-1}I_{\{X_{t_{k}},i\}}h$. 

The following theorem shows the weak consistency of the discrete type QVEs $\tilde{\sigma}_{i,n}$.

\hfill

\begin{thm}\label{thmsigma2}
	Assume that $0<\Lambda<\infty$ and $h=\Lambda/n$. Then the discrete-type QVEs $\tilde{\sigma}_{i,n}$ admit the weak consistency, i.e.,
		\begin{equation*}
			\tilde{\sigma}^{2}_{i,n}\stackrel{P}{\longrightarrow}\sigma^{2}_{i},\quad\text{as}\quad n\rightarrow\infty.
	\end{equation*}
\end{thm}
\noindent\textbf{Proof of Theorem \ref{thmsigma2}}. Let $Q^{(i)}_{\Lambda}=\{k: t_{k}\in(0,\Lambda), X_{t_{k}}\in[\theta_{i-1},\theta_{i})\}$. By Eq. (\ref{eqDelta}), we have
\begin{equation*}
	\begin{aligned}
		[XI_{\{X,i\}}]^{n}_{\Lambda}=&\sum_{k\in Q^{(i)}_{\Lambda}}\left(X_{t_{k+1}}-X_{t_{k}}\right)^{2}
		=\xi^{(i)}_{1,n}+\xi^{(i)}_{2,n}+2\xi^{(i)}_{3,n},\\
	\end{aligned}
\end{equation*}
where
\begin{equation*}
	\begin{aligned}
		\xi^{(i)}_{1,n}=&\sum_{k\in Q^{(i)}_{\Lambda}}\left(\int_{t_{k}}^{t_{k+1}}\sum_{j=1}^{m}(\beta_{j}-\alpha_{j}X_{t})I_{\{X_{t},j\}}\mathrm{d}t\right)^{2},\quad
		\xi^{(i)}_{2,n}=\sum_{k\in Q^{(i)}_{\Lambda}}\left(\sigma_{i}\Delta W_{t_{k+1}}\right)^{2},\\
		\xi^{(i)}_{3,n}=&\sum_{k\in Q^{(i)}_{\Lambda}}\left(\int_{t_{k}}^{t_{k+1}}\sum_{j=1}^{m}(\beta_{j}-\alpha_{j}X_{t})I_{\{X_{t},j\}}\mathrm{d}t\right)\left(\sigma_{i}\Delta W_{t_{k+1}}\right).
	\end{aligned}
\end{equation*}
Recall the definition of $L_{i,n}$ and $L_{i,\Lambda}$, we have
\begin{equation*}
	\lim_{n\rightarrow\infty}L_{i,n}=L_{i,\Lambda}=\int_{0}^{\Lambda}I_{\{X_{t},i\}}\mathrm{d}t,\quad a.s.
\end{equation*}
By some trivial calculations, we have
\begin{equation}\label{eqdissigma1}
	\lim_{n\rightarrow\infty}\frac{\xi_{1,n}^{(i)}}{L_{i,n}}=\mathcal{O}(h), \quad a.s.
\end{equation}
Since the quadratic variation of a martingale can be approximated by the sum of squared increments over shrinking partitions, we have
\begin{equation*}
	\begin{aligned}
		\xi^{(i)}_{2,n}&=\sum_{k=0}^{n-1}\left(\sigma_{i}I_{\{X_{t},i\}}\Delta W_{t_{k+1}}\right)^{2}\stackrel{P}{\longrightarrow}\braket{\xi_{2}^{(i)}}_{\Lambda}=\sigma_{i}^{2}\int_{0}^{\Lambda}I_{\{X_{t},i\}}\mathrm{d}t.
	\end{aligned}
\end{equation*}
Then we have
\begin{equation}\label{eqdissigma2}
	\frac{\xi_{2,n}^{(i)}}{L_{i,n}}\stackrel{P}{\longrightarrow}\sigma_{i}^{2},\quad\text{as}\quad n\rightarrow\infty.
\end{equation}
By the scaling properties of the increments of the Brownian motion, we have
\begin{equation}\label{eqdissigma3}
	\begin{aligned}
		\frac{\xi_{3,n}^{(i)}}{L_{i,n}}\stackrel{P}{\longrightarrow}0,\quad\text{as}\quad n\rightarrow\infty.
	\end{aligned}
\end{equation}
Combining Eq. (\ref{eqdissigma1})-(\ref{eqdissigma3}) yields 

	\begin{equation*}
		\frac{[XI_{\{X,i\}}]_{\Lambda}^{n}}{L_{i,n}}\stackrel{P}{\longrightarrow}\sigma_{i}^{2},\quad\text{as}\quad n\rightarrow\infty.
\end{equation*}
Thus, we complete the proof.
\hfill$\square$

\hfill

Now we consider MQVEs for the diffusion parameters $\sigma_{i}$. The main idea is to obtain the consistent estimators for the drift parameters $\alpha_{i}$ and $\beta_{i}$ first. Then we could minimize the error made by the drift term as much as possible. Hence, we need the long-time observations.
	
	When the continuously observed process $\{X_{t}\}_{0\leq t\leq T}$ is obtained, the continuous-type MQVEs for $\sigma_{i}$ are defined as 
	{\begin{equation*}
		\hat{\sigma}_{i,T}=\sqrt{\frac{\braket{Y^{(i)}}_{T}}{L_{i,T}}},
	\end{equation*}}
	where $Y^{(i)}_{t}=\int_{0}^{t}I_{\{X_{s},i\}}\mathrm{d}X_{s}-\hat{\beta}_{i,T}\int_{0}^{t}I_{\{X_{s},i\}}\mathrm{d}s+\hat{\alpha}_{i,T}\int_{0}^{t}X_{t}I_{\{X_{s},i\}}\mathrm{d}s$. 
	
	The following theorem shows the strong consistency of the continuous-type MQVEs $\tilde{\sigma}_{i,T}$.
	
	\hfill
	
	\begin{thm}\label{thmsigma4}
		{The continuous-type MQVEs $\hat{\sigma}_{i,T}$ admit the strong consistency, i.e.,
		\begin{equation*}
			\lim_{T\rightarrow\infty}\hat{\sigma}^{2}_{i,T}=\sigma^{2}_{i},\quad a.s.
		\end{equation*}}
	\end{thm}

\hfill

	\noindent\textbf{Proof of Theorem \ref{thmsigma4}}. According to (1) of Theorem \ref{thmc2}, we have 
	\begin{equation*}
		\left(\hat{\beta}_{i,T},\hat{\alpha}_{i,T}\right)\stackrel{a.s.}{\longrightarrow}\left(\beta_{i},\alpha_{i}\right),\quad\text{as}\quad T\rightarrow\infty.
	\end{equation*}
	By Eq. (\ref{eq1}), we have
	\begin{equation*}
		\int_{0}^{t}I_{\{X_{s},i\}}\mathrm{d}X_{s}=\beta_{i}\int_{0}^{t}I_{\{X_{s},i\}}\mathrm{d}s+\alpha_{i}\int_{0}^{t}X_{s}I_{\{X_{s},i\}}\mathrm{d}s+\sigma_{i}\int_{0}^{t}I_{\{X_{s},i\}}\mathrm{d}W_{s}.
	\end{equation*}
	Hence, we have
	\begin{equation*}
		Y_{t}^{(i)}=\sigma_{i}\int_{0}^{t}I_{\{X_{s},i\}}\mathrm{d}W_{s}.
	\end{equation*}
	Furthermore, the quadratic variation of $Y^{(i)}$ is given by
	\begin{equation*}
		\braket{Y^{(i)}}_{T}=\sigma^{2}_{i}\int_{0}^{T}I_{\{X_{t},i\}}\mathrm{d}t,
	\end{equation*}
	which yields
	\begin{equation*}
		\lim_{T\rightarrow\infty}\sqrt{\frac{\braket{Y^{(i)}}_{T}}{L_{i,T}}}=\sigma_{i},\quad a.s.
	\end{equation*}
	\hfill$\square$

\hfill

When the long time discretely observed process $\{X_{t_{k}}\}_{k=0}^{n}$ is obtained, i.e., $h\rightarrow0$, $n\rightarrow\infty$, and $nh\rightarrow\infty$, the discrete-type MQVEs for $\sigma_{i}$ are defined as 
\begin{equation*}
	\begin{aligned}
		&\hat{\sigma}_{i,n}
		=\sqrt{\frac{\sum_{k=0}^{n-1}\left(X_{t_{k+1}}-X_{t_{k}}-\sum_{j=1}^{m}\left(\hat{\beta}_{j,n}-\hat{\alpha}_{j,n}X_{t_{k}}\right)I_{\{X_{t_{k}},j\}}h\right)^{2}I_{\{X_{t_{k}},i\}}}{L_{i,n}}}.
	\end{aligned}
\end{equation*}

{In the following theorem, we demonstrate the consistency of the MQVE.}

\hfill

\begin{thm}\label{thmsigma3}
	Under the Hypothesis \ref{hypo1}, the discrete-type MQVEs $\hat{\sigma}_{i,n}$ admit the weak convergence, i.e.,
 {\begin{equation*}
		\hat{\sigma}^{2}_{i,n}\stackrel{P}{\longrightarrow}\sigma^{2}_{i},\quad\text{as}\quad n\rightarrow\infty.
	\end{equation*}}
\end{thm}
\noindent\textbf{Proof of Theorem \ref{thmsigma3}}. Let $Q^{(i)}_{nh}=\{k: t_{k}\in(0,nh), X_{t_{k}}\in[\theta_{i-1},\theta_{i})\}$. By Eq. (\ref{eqDelta}), we have
\begin{equation*}
	\begin{aligned}
		&\sum_{k=0}^{n-1}\left(X_{t_{k+1}}-X_{t_{k}}-\sum_{j=1}^{m}\left(\hat{\beta}_{j,n}-\hat{\alpha}_{j,n}X_{t_{k}}\right)I_{\{X_{t_{k}},j\}}h\right)^{2}I_{\{X_{t_{k}},i\}}\\
		=&\sum_{k\in Q^{(i)}_{nh}}\left(X_{t_{k+1}}-X_{t_{k}}-(\hat{\beta}_{i,n}-\hat{\alpha}_{i,n}X_{t_{k}})h\right)^{2}\\
		=&\xi^{(i)}_{1,n}+\xi^{(i)}_{2,n}+2\xi^{(i)}_{3,n},\\
	\end{aligned}
\end{equation*}
where
	\begin{align*}
		\xi^{(i)}_{1,n}=&\sum_{k\in Q^{(i)}_{nh}}\left(\int_{t_{k}}^{t_{k+1}}\sum_{j=1}^{m}(\beta_{j}-\alpha_{j}X_{t})I_{\{X_{t},j\}}-(\hat{\beta}_{i,n}-\hat{\alpha}_{i,n}X_{t_{k}})\mathrm{d}t\right)^{2},\quad\\
		\xi^{(i)}_{2,n}=&\sum_{k\in Q^{(i)}_{nh}}\left(\sigma_{i}\Delta W_{t_{k+1}}\right)^{2},\\
		\xi^{(i)}_{3,n}=&\sum_{k\in Q^{(i)}_{nh}}\left(\int_{t_{k}}^{t_{k+1}}\sum_{j=1}^{m}(\beta_{j}-\alpha_{j}X_{t})I_{\{X_{t},j\}}-(\hat{\beta}_{i,n}-\hat{\alpha}_{i,n}X_{t_{k}})h\right)\left(\sigma_{i}\Delta W_{t_{k+1}}\right).
	\end{align*}
On the one hand, similar arguments as Eq. (\ref{eqdissigma2}) and (\ref{eqdissigma3}) yield
\begin{equation*}
	\frac{\xi_{2,n}^{(i)}}{L_{i,n}}\stackrel{P}{\longrightarrow}\sigma_{i}^{2}
	\quad\text{and}\quad\frac{\xi_{3,n}^{(i)}}{L_{i,n}}\stackrel{P}{\longrightarrow}0,\quad\text{as}\quad n\rightarrow\infty.
\end{equation*}
On the other hand, we now consider the error term $\xi_{1,n}^{(i)}/L_{i,n}$, which is the main improvement of MQVEs $\hat{\sigma}_{i,n}$ compared with QVEs $\tilde{\sigma}_{i,n}$.
{Decomposing $I_{\{X_{t},i\}}$ into three terms as Eq. (\ref{eqI123}) yields
\begin{equation*}
\begin{aligned}
\xi_{1,n}^{(i)}=&\sum_{k=0}^{n-1}\left(I_{\{X_{t_{k}},i\}}^{(1)}+I_{\{X_{t_{k}},i\}}^{(2)}+I_{\{X_{t_{k}},i\}}^{(3)}\right)\left[\int_{t_{k}}^{t_{k+1}}\left(\sum_{j=1}^{m}\beta_{j}I_{\{X_{t},j\}}-\hat{\beta}_{i,n}\right)\mathrm{d}t\right.\\
&\left.+\int_{t_{k}}^{t_{k+1}}\left(\sum_{j=1}^{m}\alpha_{j}X_{t}I_{\{X_{t},j\}}-\hat{\alpha}_{i,n}X_{t_{k}}\right)\mathrm{d}t\right]^{2}.
\end{aligned}
\end{equation*}
When $I_{\{X_{t_{k}},i\}}^{(2)}=1$, by Eq. (\ref{eqd2thm1}), and (\ref{eqalphao})-(\ref{eqbetao}), we have
\begin{equation*}
\begin{aligned}
&\left[\int_{t_{k}}^{t_{k+1}}\left(\sum_{j=1}^{m}\beta_{j}I_{\{X_{t},j\}}-\hat{\beta}_{i,n}\right)\mathrm{d}t+\int_{t_{k}}^{t_{k+1}}\left(\sum_{j=1}^{m}\alpha_{j}X_{t}I_{\{X_{t},j\}}-\hat{\alpha}_{i,n}X_{t_{k}}\right)\mathrm{d}t\right]^{2}\\
\leq&\left[\int_{t_{k}}^{t_{k+1}}\left\lvert\beta_{i}-\hat{\beta}_{i,n}\right\rvert\mathrm{d}t+\int_{t_{k}}^{t_{k+1}}\left\lvert\alpha_{i}X_{t}-\alpha_{i}X_{t_{k}}+\alpha_{i}X_{t_{k}}-\hat{\alpha}_{i,n}X_{t_{k}}\right\rvert\mathrm{d}t\right]^{2}\\
\leq&\left[\left\lvert\beta_{i}-\hat{\beta}_{i,n}\right\rvert h+\int_{t_{k}}^{t_{k+1}}\left(\lvert\alpha_{i}\rvert\lvert X_{t}-X_{t_{k}}\rvert+\lvert\alpha_{i}-\hat{\alpha}_{i,n}\rvert \lvert X_{t_{k}}\rvert\right)\mathrm{d}t\right]^{2}\\
\leq&\mathcal{O}\left(\frac{h\log\log nh}{n}\right),\quad a.s.
\end{aligned}
\end{equation*}
We thus have
\begin{equation*}
\begin{aligned}
&\lim_{n\rightarrow\infty}\frac{1}{L_{i,n}}\sum_{k=0}^{n-1}I_{\{X_{t_{k}},i\}}^{(2)}\left[\int_{t_{k}}^{t_{k+1}}\left(\sum_{j=1}^{m}\beta_{j}I_{\{X_{t},j\}}-\hat{\beta}_{i,n}\right)\mathrm{d}t\right.\\
&\left.\qquad\qquad\qquad\qquad\qquad\qquad+\int_{t_{k}}^{t_{k+1}}\left(\sum_{j=1}^{m}\alpha_{j}X_{t}I_{\{X_{t},j\}}-\hat{\alpha}_{i,n}X_{t_{k}}\right)\mathrm{d}t\right]^{2}\\
\leq& \mathcal{O}\left(\frac{\log\log nh}{n}\right),\quad a.s.
\end{aligned}
\end{equation*}
By Eq. (\ref{eqergodic}), we have
	\begin{align*}
		&\lim_{n\rightarrow\infty}\frac{1}{L_{i,n}}\sum_{k=0}^{n-1}\left(I_{\{X_{t_{k}},i\}}^{(1)}+I_{\{X_{t_{k}},i\}}^{(3)}\right)\left[\int_{t_{k}}^{t_{k+1}}\left(\sum_{j=1}^{m}\beta_{j}I_{\{X_{t},j\}}-\hat{\beta}_{i,n}\right)\mathrm{d}t\right.\\
		&\left.\qquad\qquad\qquad\qquad+\int_{t_{k}}^{t_{k+1}}\left(\sum_{j=1}^{m}\alpha_{j}X_{t}I_{\{X_{t},j\}}-\hat{\alpha}_{i,n}X_{t_{k}}\right)\mathrm{d}t\right]^{2}\\
		\leq& \lim_{n\rightarrow\infty}\frac{2}{L_{i,n}}\sum_{k=0}^{n-1}\left(I_{\{X_{t_{k}},i\}}^{(1)}+I_{\{X_{t_{k}},i\}}^{(3)}\right)\left[\beta^{\ast}h+\alpha^{\ast}\sup_{s\in(t_{k}, t_{k+1}]}|X_{s}|h\right]^{2}\\
		=&\mathcal{O}\left(\sqrt{h^{3}\log\log(1/h)}\right),\quad a.s.
	\end{align*}
Hence, we have
\begin{equation*}
\lim_{n\rightarrow\infty}\frac{\xi_{i,n}^{(1)}}{L_{i,n}}\leq\mathcal{O}\left(\frac{\log\log nh}{n}\right)+\mathcal{O}\left(\sqrt{h^{3}\log\log(1/h)}\right),\quad a.s.
\end{equation*}
Hence, we complete the proof.
}
\hfill$\square$

\hfill

\begin{rmk}
	The main improvement of MQVEs compared with QVEs is to give the consistent estimators for the drift parameters, which ensures us to eliminate the error term $\xi^{(i)}_{1,n}$.
\end{rmk}

\hfill


\section{Numerical Results and Applications}\label{secnr}
\subsection{Numerical results}
In this section, we will conduct simulation studies to check the performance of the proposed estimators and compare our estimators with the GMEs proposed in \cite{Hu and Xi(2022)}. Note that it is significantly more complicated to explore the multiple thresholds case by the generalized moment approach. Hence, we consider the following three scenarios.
\begin{enumerate}[label=Scenario \arabic*,align=left]
	\item\label{S1}(Multi-regime). Set $m=3$, $\alpha=(1,2,3)$,  $\beta=(0.3,0.5,0.7)$,  $\sigma=(1,2,3)$, the two thresholds $\theta_{1}=-0.5$ and $\theta_{2}=0.5$. We study the LSEs for $\alpha$ and $\beta$, and compare the QVEs with MQVEs for $\sigma$.
	\item\label{S2}(Two-regime) Set $m=2$, $\alpha=(1,2)$,  $\beta=(0,0)$, $\theta_{1}=0$, and $\sigma=1$. We compare our estimators with the GMEs for $\alpha$.
	\item\label{S3}(Two-regime) Set $m=2$, $\alpha=(1,2)$, $\beta=(-0.3,0.3)$, $\theta_{1}=0$, and $\sigma=1$. We compare our estimators with the GMEs for $\alpha$ and $\beta$.
\end{enumerate}

\begin{table}[h!]
	\caption{Mean and standard deviation of $\hat{\alpha}_{i,n}$ and $\hat{\beta}_{i,n}$ for $i = 1, 2, 3$ with $h = 0.1$ under different choices of sample sizes $n$.}
	\begin{tabular}{ccccccc}		
		\hline
		&&\multicolumn{5}{l}{Sample sizes $n$}
		\tabularnewline
		\hline
		\makebox[0.05\textwidth]{}&
		\makebox[0.10\textwidth]{}&
		\makebox[0.12\textwidth]{1000}&  
		\makebox[0.12\textwidth]{2000}&
		\makebox[0.12\textwidth]{3000}&  
		\makebox[0.12\textwidth]{4000}&
		\makebox[0.12\textwidth]{5000}
		\tabularnewline
		\hline
		$\alpha_{1}$ 
		& bias    & 0.082 & 0.020 & 0.021 & 0.004 & 0.007 \tabularnewline
		& Std.dev & 0.412 & 0.291 & 0.247 & 0.198 & 0.183 \tabularnewline
		$\alpha_{2}$ 
		& bias    & 0.054 & -0.026 & 0.030 & 0.032 & 0.018 \tabularnewline
		& Std.dev & 1.270 & 0.878 & 0.693 & 0.626 & 0.549 \tabularnewline
		$\alpha_{3}$ 
		& bias    & 0.150 & 0.079 & 0.039 & 0.020 & 0.023 \tabularnewline
		& Std.dev & 0.917 & 0.625 & 0.507 & 0.439 & 0.377 \tabularnewline
		$\beta_{1}$ 
		& bias    & -0.064 & -0.011 & -0.013 & -0.001 & -0.004 \tabularnewline
		& Std.dev & 0.432 & 0.308 & 0.262 & 0.215 & 0.194 \tabularnewline
		$\beta_{2}$ 
		& bias    & -0.001 & -0.001 & 0.008 & -0.008 & -0.008 \tabularnewline
		& Std.dev & 0.372 & 0.250 & 0.215 & 0.177 & 0.169 \tabularnewline
		$\beta_{3}$ 
		& bias    & 0.109 & 0.088 & 0.024 & 0.007 & 0.019 \tabularnewline
		& Std.dev & 1.344 & 0.954 & 0.785 & 0.671 & 0.576 \tabularnewline
		\hline 
	\end{tabular}
	\label{table1}
\end{table}

\begin{figure}[h!]
	\centering
	\subfigure[Normal QQ-plot for  $\hat{\alpha}_{1}$]{
		\includegraphics[scale=0.26]{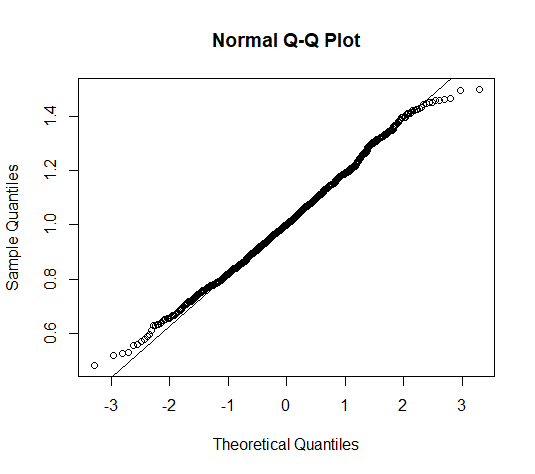}
	}
	\subfigure[Normal QQ-plot for  $\hat{\alpha}_{2}$]{
		\includegraphics[scale=0.26]{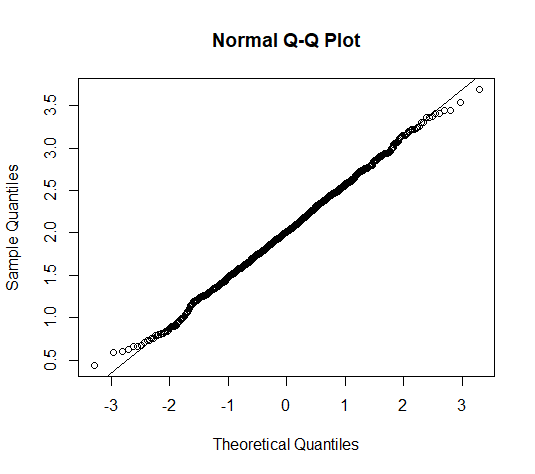}
	}
	\subfigure[Normal QQ-plot for  $\hat{\alpha}_{3}$]{
		\includegraphics[scale=0.26]{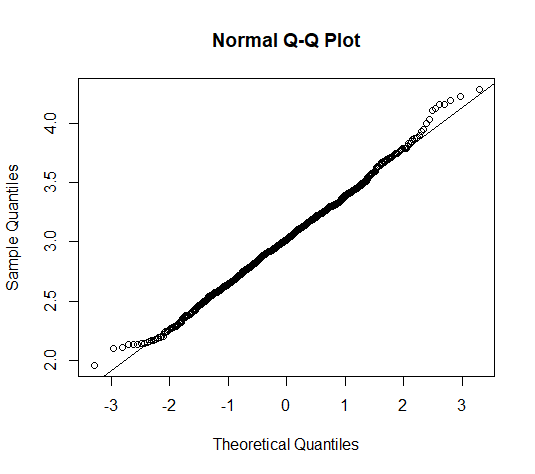}
	}
	\subfigure[Normal QQ-plot for  $\hat{\beta}_{1}$]{
		\includegraphics[scale=0.26]{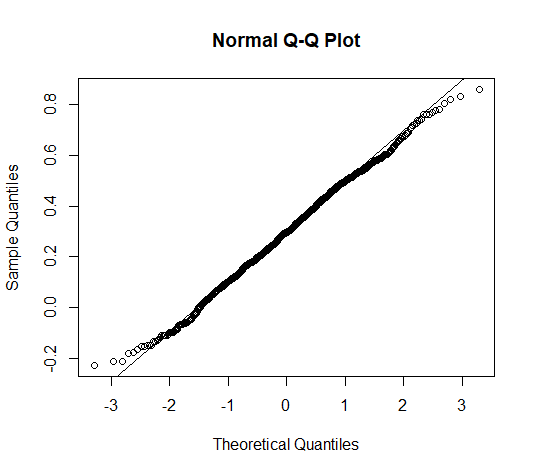}
	}
	\subfigure[Normal QQ-plot for  $\hat{\beta}_{2}$]{
		\includegraphics[scale=0.26]{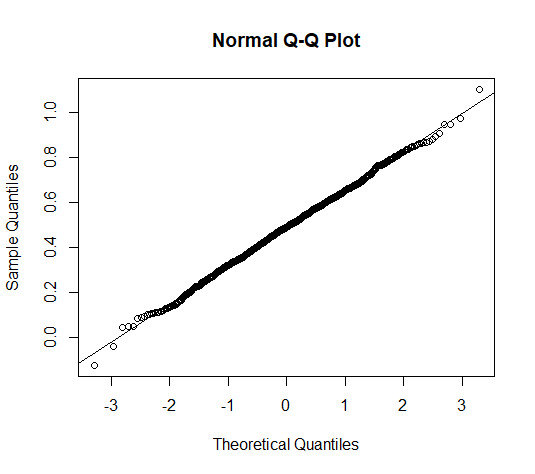}
	}
	\subfigure[Normal QQ-plot for  $\hat{\beta}_{3}$]{
		\includegraphics[scale=0.26]{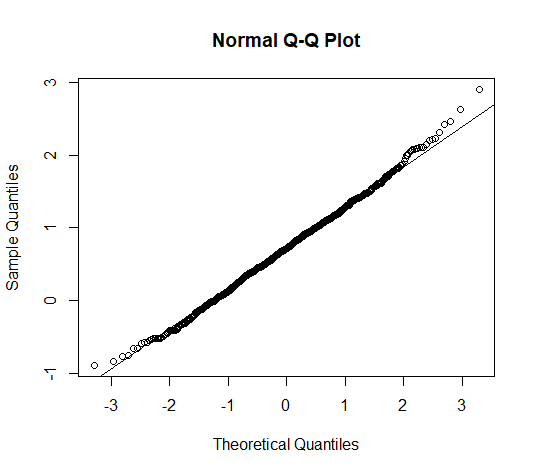}
	}
	\caption{Normal QQ plot for $1000$ samples of the drift parameters with $n = 5000$ and $h=0.1$.}
	\label{figqq}
\end{figure}

\begin{table}[h!]
	\caption{Mean and standard deviation of  $\tilde{\sigma}_{i,n}$ and $\hat{\sigma}_{i,n}$ for $i = 1, 2, 3$  with $h=0.01$ and $h = 0.1$ under different choices of sample sizes $n$, where $\sigma$ denotes the results of MQVEs and $\sigma^{\prime}$ denotes the results of QVEs.}
	\begin{tabular}{cccccccc}		
		\hline
		&&&\multicolumn{5}{l}{Sample sizes $n$}
		\tabularnewline
		\hline
		\makebox[0.045\textwidth]{}&
		\makebox[0.045\textwidth]{}&
		\makebox[0.09\textwidth]{}&
		\makebox[0.09\textwidth]{1000}&  
		\makebox[0.09\textwidth]{2000}&
		\makebox[0.09\textwidth]{3000}&  
		\makebox[0.09\textwidth]{4000}&
		\makebox[0.09\textwidth]{5000}
		\tabularnewline
		\hline
		$h=0.01$ & $\sigma_{1}$ 
		& bias    & -0.016 & -0.006 & -0.005 & -0.003 & -0.001 \tabularnewline
		&& Std.dev & 0.049 & 0.028 & 0.021 & 0.018 & 0.015 \tabularnewline
		&$\sigma^{\prime}_{1}$ 
		& bias    & -0.008 & -0.008 & -0.009 & -0.008 & -0.007 \tabularnewline
		&& Std.dev & 0.036 & 0.026 & 0.020 & 0.018 & 0.015 \tabularnewline
		&$\sigma_{2}$ 
		& bias    & 0.003 & 0.002 & -0.001 & 0.000 & -0.001 \tabularnewline
		&& Std.dev & 0.079 & 0.054 & 0.046 & 0.038 & 0.035 \tabularnewline
		&$\sigma^{\prime}_{2}$ 
		& bias    & -0.001 & -0.001 & -0.003 & -0.002 & -0.004 \tabularnewline
		&& Std.dev & 0.079 & 0.054 & 0.046 &  0.038 & 0.035 \tabularnewline
		&$\sigma_{3}$ 
		& bias    & -0.045 & -0.017 & -0.008 & -0.005 & -0.003  \tabularnewline
		&& Std.dev & 0.207 & 0.111 & 0.091 & 0.075 & 0.065 \tabularnewline
		&$\sigma^{\prime}_{3}$ 
		& bias    & -0.028 & -0.033 & -0.027 & -0.027 & -0.026 \tabularnewline
		&& Std.dev & 0.166 & 0.109 & 0.090 & 0.074 & 0.065 \tabularnewline
		\hline 
		$h=0.1$ &$\sigma_{1}$ 
		& bias    & 0.006 & 0.003 & 0.002 & 0.001 & 0.001 \tabularnewline
		&& Std.dev & 0.035 & 0.026 & 0.019 & 0.017 & 0.015 \tabularnewline
		&$\sigma^{\prime}_{1}$ 
		& bias    & 0.089 & 0.088 & 0.089 & 0.088 & 0.088 \tabularnewline
		&& Std.dev & 0.035 & 0.026 & 0.020 & 0.018 & 0.016 \tabularnewline
		&$\sigma_{2}$ 
		& bias    & -0.005 & -0.002 & 0.001 & -0.001 & 0.001 \tabularnewline
		&& Std.dev & 0.078 & 0.057 & 0.045 & 0.038 & 0.034 \tabularnewline
		&$\sigma^{\prime}_{2}$ 
		& bias    & 0.016 & 0.017 & 0.021 & 0.018 & 0.021 \tabularnewline
		&& Std.dev & 0.079 & 0.058 & 0.046 & 0.039 & 0.035 \tabularnewline
		&$\sigma_{3}$ 
		& bias    & 0.008 & 0.007 & 0.006 & 0.002 & 0.002 \tabularnewline
		&& Std.dev & 0.141 & 0.098 & 0.081 & 0.072 & 0.063 \tabularnewline
		&$\sigma^{\prime}_{3}$ 
		& bias    & 0.260 & 0.264 & 0.262 & 0.260 & 0.261 \tabularnewline
		&& Std.dev & 0.150 & 0.104 & 0.088 & 0.078 & 0.068 \tabularnewline
		\hline
	\end{tabular}
	\label{tablesigma}
\end{table}

\begin{table}[h!]
	\caption{Summary statistics comparing our method and the generalized moment method: mean and standard deviation of $\hat{\alpha}_{i,n}$ for $i = 1, 2$ with $h = 0.1$ under different choices of sample sizes $n$.}
	\begin{tabular}{cccccc}
		\hline
		&&\multicolumn{2}{c}{Least squares estimator}
		&\multicolumn{2}{c}{Generalized moment estimator} 
		\tabularnewline
		\hline
		\makebox[0.04\textwidth]{n}&
		\makebox[0.04\textwidth]{}&
		\makebox[0.17\textwidth]{bias}&  
		\makebox[0.17\textwidth]{Std.dev}&
		\makebox[0.17\textwidth]{bias}&  
		\makebox[0.17\textwidth]{Std.dev}
		\tabularnewline
		1000 &$\alpha_{1}$ & 0.034 & 0.186 & 0.006 & 0.203  \tabularnewline
		&$\alpha_{2}$ & 0.061 & 0.321 & -0.123 & 0.375\tabularnewline
		2000 &$\alpha_{1}$ & 0.012 & 0.131 & -0.030 & 0.138 \tabularnewline
		&$\alpha_{2}$ & 0.022 & 0.211 & -0.141 & 0.248\tabularnewline
		3000 &$\alpha_{1}$ & 0.012 & 0.104 & -0.034 & 0.111 \tabularnewline
		&$\alpha_{2}$ & 0.013 & 0.177 & -0.165 & 0.210 \tabularnewline
		4000 &$\alpha_{1}$ & 0.018 & 0.090 & -0.045 & 0.096 \tabularnewline
		&$\alpha_{2}$ & 0.016 & 0.150 & -0.172 & 0.175 \tabularnewline
		5000  &$\alpha_{1}$ & 0.018 & 0.079 & -0.051 & 0.083 \tabularnewline
		&$\alpha_{2}$ & 0.018 & 0.138 & -0.175 & 0.158 \tabularnewline
		\hline 
	\end{tabular}
	\label{table2}
\end{table}

\begin{table}[h!]
	\caption{Summary statistics comparing our method and the generalized moment method: mean and standard deviation of $\hat{\alpha}_{i,n}$ and $\hat{\beta}_{i,n}$ for $i = 1, 2$ with $h = 0.1$ under different choices of sample sizes $n$.}
	\begin{tabular}{cccccc}
		\hline 
		&&\multicolumn{2}{c}{Least squares estimator}
		&\multicolumn{2}{c}{Generalized moment estimator} 
		\tabularnewline
		\hline
		\makebox[0.04\textwidth]{n}&
		\makebox[0.04\textwidth]{}&
		\makebox[0.17\textwidth]{bias}&  
		\makebox[0.17\textwidth]{Std.dev}&
		\makebox[0.17\textwidth]{bias}&  
		\makebox[0.17\textwidth]{Std.dev}
		\tabularnewline
		1000 &$\alpha_{1}$ & 0.031 & 0.197 & -0.345 & 0.136   \tabularnewline
		&$\alpha_{2}$ & 0.101 & 0.478 & -0.056 & 0.562 \tabularnewline
		&$\beta_{1}$ & -0.006 & 0.075 & 0.337 & 0.084 \tabularnewline
		&$\beta_{2}$ & 0.008 & 0.110 & -0.030 & 0.293 
		\tabularnewline
		2000 &$\alpha_{1}$ & 0.020 & 0.140 & -0.349 & 0.093 \tabularnewline
		&$\alpha_{2}$ & 0.061 & 0.336 & -0.127 & 0.406 \tabularnewline
		&$\beta_{1}$ & -0.005 & 0.055 & 0.320 & 0.044 \tabularnewline
		&$\beta_{2}$ & 0.005 & 0.079 & 0.030 & 0.255  \tabularnewline
		3000 &$\alpha_{1}$ &  0.018 & 0.116 & -0.349 & 0.073 \tabularnewline
		&$\alpha_{2}$ & 0.028 & 0.260 & -0.178 & 0.342 \tabularnewline
		&$\beta_{1}$ & -0.005 & 0.044 & 0.313 & 0.025 \tabularnewline
		&$\beta_{2}$ & 0.002 & 0.062 & 0.054 & 0.223  \tabularnewline
		4000 &$\alpha_{1}$ & 0.011 & 0.101 & -0.345 & 0.065 \tabularnewline
		&$\alpha_{2}$ & 0.015 & 0.226 & -0.200 & 0.293 \tabularnewline
		&$\beta_{1}$ & -0.003 & 0.039 & 0.310 & 0.013 \tabularnewline
		&$\beta_{2}$ & 0.000 & 0.055 & 0.076 & 0.189  \tabularnewline
		5000 &$\alpha_{1}$ & 0.005 & 0.086 & 0.349 & 0.057 \tabularnewline
		&$\alpha_{2}$ & 0.013 & 0.205 & 0.198 & 0.267 \tabularnewline
		&$\beta_{1}$ & -0.001 & 0.033 & 0.309 & 0.012 \tabularnewline
		&$\beta_{2}$ & 0.001 & 0.050 & 0.076 & 0.171  \tabularnewline
		\hline 
	\end{tabular}
	\label{table3}
\end{table}

We use the Euler scheme to simulate the threshold OU process by discretizing Eq. (\ref{eq1}) on the interval $[0, nh]$ with different mesh sizes $h$ and sample sizes $n$.
For each scenario, we generate $N=1000$ Monte Carlo simulations of sample paths and each path consists of $n=1000, 2000, 3000, 4000$, and $5000$ observations. The overall parameter estimates are evaluated by the bias and standard deviation (Std.dev). 

Table \ref{table1} summarizes the main findings of \ref{S1} over $1000$ Monte Carlo simulations. We observe that as the sample size increases, the bias decreases and is small. The empirical and model-based standard deviation agree reasonably well. The performance of the LSEs $\alpha$ and $\beta$ improves with larger sample sizes in the multi-regime cases. 

Table \ref{tablesigma} summarizes the results of the QVEs and MQVEs. We can see from the table that the QVEs and MQVEs all have good performance when the mesh size is small($h=0.01$). However, when the mesh size is large($h=0.1$), the QVEs may lead to a rather large bias and such bias does not vanish as the sample size increases. Hence, we conclude that MQVEs have improved performance relative to the QVEs, especially when the sample size and mesh size are large.

The normal QQ-plots for the estimators of 1000 Monte Carlo simulations with $n=5000$ are presented in Figure \ref{figqq}. We can see from the figure that the estimators admit the asymptotic normality, which supports the results in Theorem \ref{thmc1}-\ref{thmd2}.

Shown in Table \ref{table2} and \ref{table3} are the mean and standard deviation of the LSEs and GMEs in \ref{S2} and \ref{S3}. The results exhibit that a rather large bias may be incurred by the generalized moment method. As the sample size increases, such biases may not attenuate. {When employing the least squares method, the bias is significantly reduced compared to the generalized moment method.} Heuristically, the two estimators perform better as the sample size becomes larger, but for relatively small sample sizes, the LSEs outperform the GMEs.

	\begin{figure}[h!]
		\centering
		\subfigure[10-year treasury rate]{
			\includegraphics[scale=0.3]{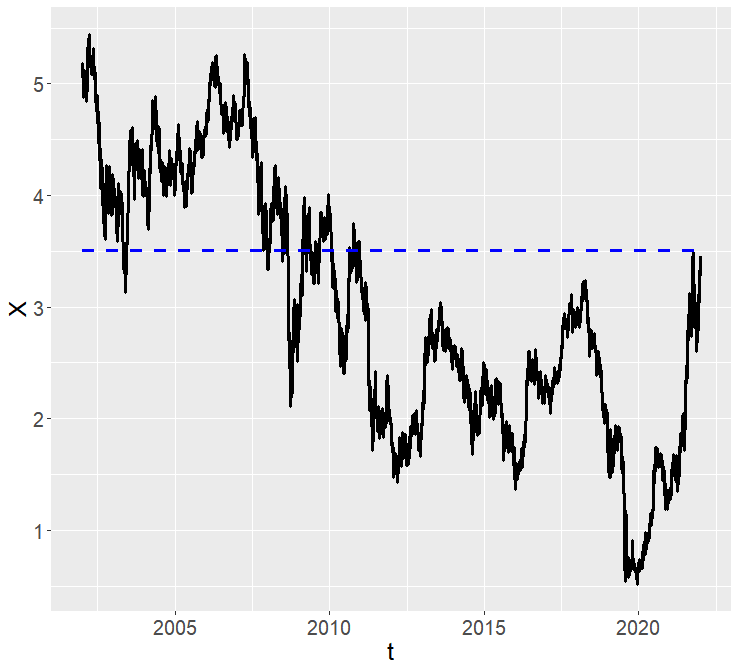}
		}
		\subfigure[5-year treasury rate]{
			\includegraphics[scale=0.3]{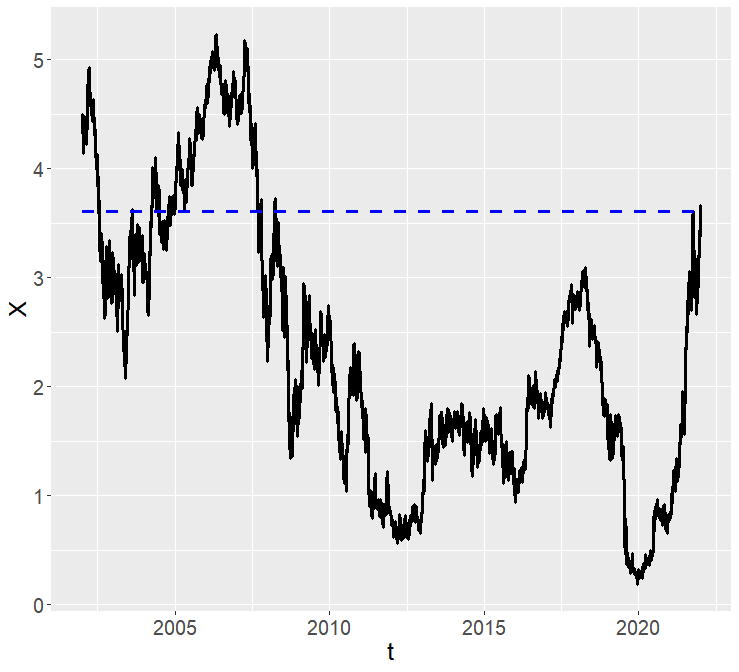}
		}
		\subfigure[2-year treasury rate]{
			\includegraphics[scale=0.3]{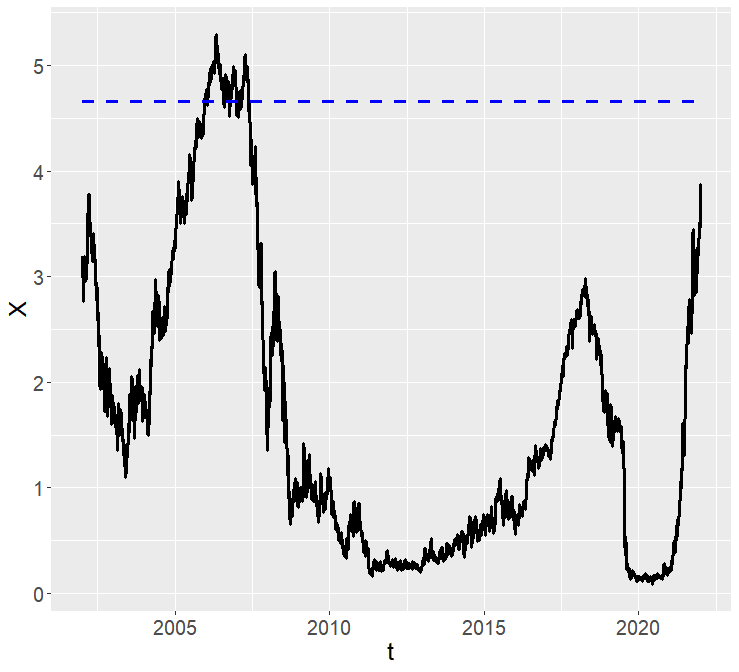}
		}
		\subfigure[1-year treasury rate]{
			\includegraphics[scale=0.3]{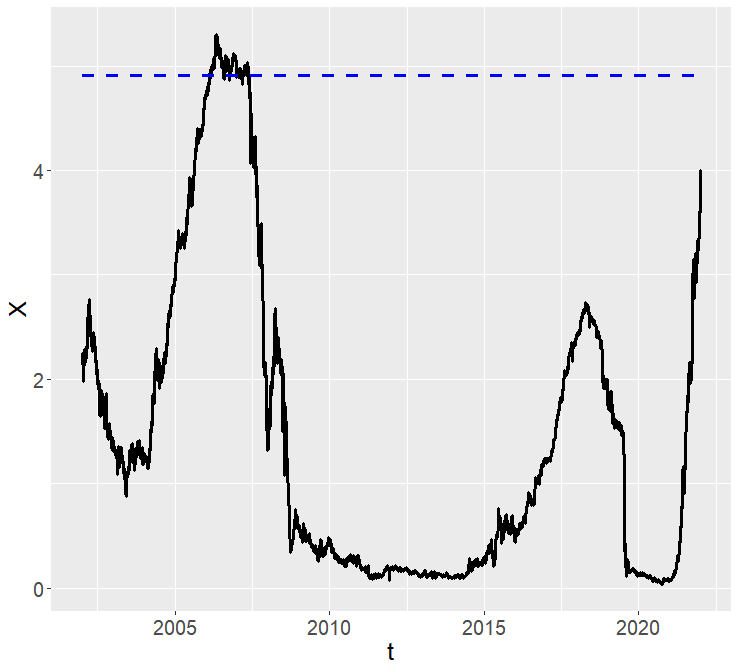}
		}
		\caption{Daily U.S. treasury rate (solid black line) with a horizontal dashed blue line being the thresholds.}
		\label{fig1}
	\end{figure}
	
	\begin{table}[h!]
		\caption{Estimating results of the U.S. treasury rate.}
		\begin{tabular}{ccccc}
			\hline 
			&\multicolumn{2}{c}{10-year treasury rate}
			&\multicolumn{2}{c}{5-year treasury rate} 
			\tabularnewline
			\hline
			\makebox[0.065\textwidth]{}
			&\makebox[0.1\textwidth]{value}
			&\makebox[0.28\textwidth]{CI}  
			&\makebox[0.1\textwidth]{value}
			&\makebox[0.28\textwidth]{CI} 
			\tabularnewline
			\hline 
			$\theta$ & 3.510 &  & 3.610 &  \tabularnewline
			$\alpha_{1}$ & 0.018 & [0.009, 0.027] & 0.036 & [0.024, 0.048]  \tabularnewline
			$\beta_{1}$ & 0.058 & [-0.002, 0.119] & 0.075 & [0.018, 0.132]  \tabularnewline
			$\alpha_{2}$ & 0.062 & [0.058, 0.067] & 0.006 & [0.003, 0.009]  \tabularnewline
			$\beta_{2}$ & 0.215 & [0.180, 0.251] & -0.023 & [-0.055, 0.008]  \tabularnewline
			$\sigma_{1}$ & 0.254 &  & 0.279 &   \tabularnewline
			$\sigma_{2}$ & 0.280 &  & 0.250 &   \tabularnewline
			\hline
			&\multicolumn{2}{c}{2-year treasury rate}
			&\multicolumn{2}{c}{1-year treasury rate}
			\tabularnewline
			\hline 
			\makebox[0.065\textwidth]{}
			&\makebox[0.1\textwidth]{value}
			&\makebox[0.28\textwidth]{CI}  
			&\makebox[0.1\textwidth]{value}
			&\makebox[0.28\textwidth]{CI} 
			\tabularnewline
			\hline
			$\theta$ & 4.660 &  & 4.907 &  \tabularnewline
			$\alpha_{1}$ & 0.006 & [-0.006, 0.019] & -0.005 & [-0.015, 0.006] \tabularnewline
			$\beta_{1}$ & 0.016 & [-0.029, 0.061] & 0.004 & [-0.028, 0.036] \tabularnewline
			$\alpha_{2}$ & 0.299 & [0.298, 0.300] & 0.540 & [0.539, 0.541] \tabularnewline
			$\beta_{2}$ & 1.400 & [1.370, 1.420] & 2.678 & [2.661, 2.695] \tabularnewline
			$\sigma_{1}$ & 0.235 &  & 0.175 &  \tabularnewline
			$\sigma_{2}$ & 0.195 &  & 0.134 &  \tabularnewline
			\hline
		\end{tabular}
		\label{table4}
\end{table}

\subsection{An application in U.S. treasury rate}
In the previous discussion we presented an efficient estimation procedure for the drift parameters of a threshold OU process. \cite{Su and Chan(2015)} apply the threshold OU process to the three-month U.S. treasury rate, which is a high-frequency observed process with data collected on a daily basis. The rates are published on business days and the data shall be treated as equally spaced. However, it is not clear to verify whether the information accumulates twice over two non-business days as fast as over a one-day gap between consecutive daily data.

Here, we consider the U.S. treasury rate based on the data sets from Federal Reserve Bank. We report the treasury rate process for four assets including the ten-year, five-year, two-year, and one-year US treasury rate,  
all of which are collected from ``2002-1-2" to ``2022-9-16". We adopt the convention that the equal time interval for the ``daily" rates is $h=0.046$ while one unit in time represents one month.

	We model the rate process using the following SDE 
	\begin{equation*}
		\left\{
		\begin{aligned}
			&\mathrm{d}X_{t}=\sum_{j=1}^{m}\left[(\beta_{j}-\alpha_{j}X_{t})\mathrm{d}t+\sigma_{j}\mathrm{d}W_{t}\right]I_{\{X_{t},j\}},\\
			&X_{0}=x.
		\end{aligned}
		\right.
	\end{equation*} 
	Note that determining the number and values of thresholds (i.e., the values of $m$ and $\theta_{j}$) is an essential problem. In real-world applications, we shall transform the model as the following piecewise linear model
	\begin{equation*}
		Y_{t_{k}}=\sum_{j=1}^{m}\left[(\beta_{j}-\alpha_{j}X_{t_{k}})h+\sigma_{j}(W_{t_{k+1}}-W_{t_{k}})\right]I_{\{X_{t_{k}},j\}},
	\end{equation*}
	where $Y_{t_{k}}=X_{t_{k+1}}-X_{t_{k}}$. Then we can use the ``\emph{segmented}'' function in ``\emph{segmented}'' package of $\emph{R}$ to obtain the number and values of thresholds. {The ``\emph{segmented}'' package is a tool designed for segmented regression analysis. This package facilitates the detection and fitting of breakpoints in data, enabling the modeling of distinct linear relationships in different segments. The ``\emph{segmented}'' function within the ``\emph{segmented}'' package is primarily employed for segmented regression analysis.}
	
	Figure \ref{fig1} displays the U.S. treasury rate process and we set the thresholds as the horizontal dashed blue lines. We shall see from the figure that each of the treasury rate processes has only one threshold ($m=2$). The value of thresholds of 10-year and 5-year treasury rate processes are similar, as are 2-year and 1-year treasury rate processes. These are the ones that match the reality, since the trends of 10-year and 5-year treasury rate processes are about the same, while 2-year and 1-year treasury rate processes are about the same.

Table \ref{table4} summarizes the estimations of the drift and diffusion parameters where ``value" is the estimation and 
``CI" is the $95\%$ confidence interval. Since the long rate is closely related to economic growth and inflation, the drift term does not have a significant slope with respect to its current treasury rate. The drift term in the first regime of the one-year and two-year treasury rate process is close to $0$, which shows that the short rate evolves as a martingale process until it hits the second regime. In the second regime, the drift term has a negative slope which ensures the ergodicity of the process. In general, the short-rate process has the same general trend as the long-rate process.

\section{Disscussion}\label{seccon}
We have considered the parameter estimation problem for the parameters of a threshold Ornstein$\mathit{-}$Uhlenbeck process with multiple thresholds. Under the condition that the thresholds are known, we propose asymptotically consistent estimators for the drift and diffusion parameters. Our simulation results demonstrate that the proposed estimators are superior to the generalized moment estimators proposed by \cite{Hu and Xi(2022)} and the usual quadratic variation estimators.

There are several important directions for further research. We would like to study the least squares estimator for all parameters including the thresholds. Consider the following objective function
\begin{equation*}
	\begin{aligned}
		\min_{\alpha_{j},\theta_{j}}\int_{0}^{T}\bigg|\sum_{j=1}^{m}\left(\dot{X}_{t}+\alpha_{j}X_{t}\right)I_{\{X_{t},j\}}\bigg|^{2}\mathrm{d}t.
	\end{aligned}
\end{equation*}
By the fact that $I_{\{X_{t},j\}}=I_{\{X_{t}<\theta_{j}\}}-I_{\{X_{t}<\theta_{j-1}\}}$, the above objective function is equivalent to 
\begin{equation*}
	\begin{aligned}
		\min_{\alpha_{j},\theta_{j}}\int_{0}^{T}\bigg\lvert\sum_{j=2}^{m}(\alpha_{j-1}-\alpha_{j})X_{t}I_{\{X_{t}<\theta_{j-1}\}}+(\dot{X}_{t}+\alpha_{m}X_{t})\bigg\rvert^2\mathrm{d}t.
	\end{aligned}
\end{equation*}
It is almost impossible to obtain the closed-form formula of the estimator $\hat{\alpha}_{j}$ and $\hat{\theta}_{j}$, $j=1, \cdots, n$. However, various specialized methods have been developed to demonstrate the consistency of the suggested estimators \citep[see, e.g.,][]{Frydman(1980),van der Vaart(1998)}.
Some other further research may include investigating the statistical inference for the generalized threshold diffusions.

\backmatter

\bmhead{Acknowledgments}

The authors thank Tiefeng Jiang for the insightful comments. We also thank Yaozhong Hu and Yuejuan Xi for their enthusiastic help including the code of their method. This work is partially supported by the National Natural Science Foundation of China (No. 11871244), and the Fundamental Research Funds for the Central Universities, JLU.

\bmhead{Conflict of interest}

All authors disclosed no relevant relationships.

\bigskip

\begin{appendices}

	
	
	
\end{appendices}


\end{document}